\documentclass[10pt,a4paper]{amsart}
\usepackage{amsmath}
\usepackage{amsthm}
\usepackage{amsfonts}
\usepackage{amssymb}
\usepackage{a4,a4wide}
\usepackage{leqno}
\usepackage{color}

\usepackage{graphics}

\usepackage{color}
\usepackage[T1]{fontenc}
\usepackage{palatino}

\newtheorem{theorem}{Theorem}[section]
\newtheorem{lemma}[theorem]{Lemma}

\newtheorem{claim}[theorem]{Claim}

\theoremstyle{definition}

\theoremstyle{remark}
\newtheorem{remark}[theorem]{Remark}

\numberwithin{equation}{section}
\def\O{{\Omega}}
\def\o{{\omega}}

\def\eps{{\epsilon}}

\def\e{{\mathcal{E}}}

\def\h{{\mathcal{H}}}
\def\I{{\mathcal{I}}}

\def\n{{\mathcal{N}}}

\def\t{{\mathcal{T}}}

\def\D{{\mathcal{D}}}

\def\R{{\mathbb{R}}}

\def\N{{\mathbb{N}}}

\newcommand{\dem}[1]{\vskip 0.2\baselineskip \noindent {\bf{#1}}\vskip 0.2\baselineskip }
\newcommand{\fdem}{\vskip 0.2 pt \hfill $\square$ }

\newcommand{\transposee}[1]{{\vphantom{#1}}^{\mathit t}{#1}}
\definecolor{colorjer}{rgb}{0,0.6,0}

\def\tilde{\widetilde}

\title{Convergence to the equilibrium in a Lotka-Volterra ode competition system with mutations}
\author{J\'er\^ome Coville}
\address{\noindent J. Coville -- INRA PACA, Equipe BIOSP, Centre de Recherche d'Avignon, Domaine Saint
Paul, Site Agroparc, 84914 Avignon cedex 9, France}
\email{jerome.coville@avignon.inra.fr}

\author{Fr\'ed\'eric Fabre}
\address{\noindent F. Fabre --  INRA PACA, UR 407 Pathologie V\'eg\'etale, Centre de Recherche d'Avignon, Site Saint-Maurice,
CS 60094,
84143 Monfavet Cedex, France}
\email{frederic.fabre@avignon.inra.fr}
\date{\today}
\begin{document}


%
%
\begin{abstract}
In this paper we are investigating  the long time behaviour of the solution of a mutation competition  model of Lotka-Volterra's type.  
Our main motivation comes from the analysis of the  Lotka-Volterra's competition  system with mutation which simulates the demo-genetic dynamics of diverse virus in their host  :
$$
 \frac{dv_{i}(t)}{dt}=v_i\left[r_i-\frac{1}{K}\Psi_i(v)\right]+\sum_{j=1}^{N} \mu_{ij}(v_j-v_i). $$
  In a first part we analyse the case  where the   competition terms $\Psi_i$ are independent of the virus type $i$. In  this situation and under some rather general assumptions on the functions $\Psi_i$, the coefficients $r_i$ and the mutation matrix $\mu_{ij}$ we prove the existence of a unique positive globally stable stationary  solution i.e. the solution attracts the  trajectory  initiated  from  any nonnegative initial datum. Moreover the unique  steady state $\bar v$ is strictly positive in the sense that $\bar v_i>0$ for all $i$.       
  These results are in sharp contrast with the behaviour of Lotka-Volterra without mutation term where it is known that multiple non negative stationary solutions exist and an exclusion principle occurs (i.e For all  $i\neq i_0,  \bar v_{i}=0$ and $\bar v_{i_0}>0$ ). 
  Then  we explore a typical example that has been proposed to explain some experimental data. For such particular  models we  characterise the speed of convergence to the equilibrium. 
 In a second part, under some additional assumption, we prove the existence of a positive steady state for the full system and we analyse the long term dynamics.  The proofs  mainly rely on   the construction of a relative entropy which plays the role of a Lyapunov functional.  
  
\end{abstract}
\maketitle

\bigskip


\bigskip


\section{Introduction}

In this paper we are investigating  the long time behaviour of the solution of some   models that have been recently used in epidemiology.  Our analysis  focuses  on a Lotka-Volterra competition  system with mutation which basically simulates  the demo-genetic dynamic of a genetically diverse virus population in its hosts, highliting the numerous links existing between ecological
and within-host infection dyanmics. Such type of model  has  been  proposed to explain some experimental data e.g.\cite{Cuevas2003,Fabre2012,Miralles2001}. To be more specific the demo-genetic dynamic is modelled  by $N$ ordinary differential equations  which simulates at host scale the dynamics of  $v_i(t)$ the number of virus particles of genotype $i$ at time $t$ : 
 
 \begin{equation}
 \frac{dv_{i}(t)}{dt}=v_i\left[r_i-\frac{1}{K}\Psi_i(v)\right]+\sum_{j=1}^{N} \mu_{ij}(v_j-v_i) \label{mutselecedo-eq-sysgen}
 \end{equation}
   where $r_i, K $ and $\mu_{ij}$ represent respectively  the growth rate for each genotype, the total carrying capacity of the host and  a
   nonnegative  matrix characterising the rate of mutation from   a virus of genotype $i$ to a virus of genotype $j$.
    For each $i$, $\Psi_i(v): \R^N\to \R$ is a locally Lipschitz  application describing  the intensity of the interaction between a virus of genotype $i$ with all its competitors. 
    
 The mutation term of the system can also be interpreted as a dispersal term.  Indeed, mutation naturally corresponds to dispersal into the discrete space of genotype. But the mutation term can also handle dispersal between physical patches. With this in mind,  the above system of equation  can also be used to model  the demo-genetic dynamic of a diverse virus population in structured hosts, each host tissue types being virus "habitats" connected to each others by dispersal via fluid flow (e.g. \cite{Orive2005}). Ways of derive results for this interpretation  are discussed  in the biological  comments  subsection.
 

   
 In what follows we will always make the following assumptions on $r_i$, $ \Psi_i$ and $ \mu_{ij}$
\begin{align}
\label{mutselecedo-hyp1}
\left\{
\begin{aligned}
&\text{For all  $i$, $r_i>0$,}
\\
&\text{The matrix $(\mu_{ij})$ is  nonnegative  symmetric and irreducible }\\
&\text{ $\Psi_i(v) \in C_{loc}^{0,1}(\R^N,\R), \Psi_i(0)=0$}\\
&\text{ $\Psi_i$ is monotone increasing with respect to the natural order of $\R^N$}
\end{aligned}
\right.
\end{align} 
 
 Furthermore we will assume that for all $i$ there exists  positive constants $R_i,k_i, c_{i}$ with $k_i> 0$ so that the function   $\Psi_i$ satisfies for all  $v \in \R^{N,+}\setminus Q_{R_i}(0)$,
 \begin{align}
\label{mutselecedo-hyp2}
c_{i}\left(\sum_{j=1}^Nv_j\right)^{k_i}\le \Psi_i(v),
\end{align} 
where $Q_{R_i}(0)$ denotes the ball of radius $R_i$ and centred at $0$ associated to the $l^1$ norm.  

A typical example of model satisfying our assumption is given by 

$$  \frac{dv_{i}(t)}{dt}=v_i\left[r_i-\frac{1}{K}\sum_{j=1}^{4}\beta_{ij}v_j\right]+Mv$$

with   
$$M:= \begin{pmatrix}  (1-\mu)^2-1  &\mu(1-\mu) &\mu(1-\mu) &\mu^2 \\
     \mu(1-\mu)&(1-\mu)^2-1&  \mu^2&\mu(1-\mu) \\
    \mu(1-\mu) &\mu^2 &(1-\mu)^2-1   &\mu(1-\mu)\\
    \mu^2&\mu(1-\mu)  &\mu(1-\mu) &(1-\mu)^2 -1 \\
\end{pmatrix},
$$  
where $\mu$ is a parameter giving the point mutation rate per replication cycle and per nucleotide. This mutation matrix  corresponds to a viral population composed  of $4$ variants  differing only by one or two substitutions involved in adaptative  proprieties (e.g. pathogenicity). The interactions terms $\Psi_i:=\sum_{j=1}^{4}\beta_{ij}v_j$ can handle a wide range of possible inter-specific (inter-variants) competition rates between any pairs of virus variants.

 This particular structure  has been  used  recently to model the adaptation of plant virus to resistance genes
 \cite{Fabre2009,Fabre2012}. This particular  form of competition  is commonly used to model  virus evolution,\cite{Cuevas2003,Lafforgue2011,Miralles2001}. 
Without the mutation's matrix $(\mu_{ij})$, the system \eqref{mutselecedo-eq-sysgen} reduces to a classical competition system in the sense of Hirsch \cite{Hirsch1982,Hirsch1985, Hirsch1988}
 \begin{equation}
\frac{dv_{i}(t)}{dt}=v_i\left[r_i-\frac{1}{K}\Psi_i(v)\right]. \label{mutselecedo-eq-sysgen2}
\end{equation}  

Such  systems has been intensively studied  and many aspects are now  well understood see for example \cite{Akin1979, Burger2000,Burger1994,Champagnat2010,Crow1970,Diekmann2005,Hadeler1981,Hirsch1982,Hirsch1985, Hirsch1988, Hofbauer1985, Hofbauer1987, Jabin2011,Li1999, Perthame2007} and references therein.  In particular, the existence of stationary solution and the  asymptotic behaviour of the solution  has been  obtained in \cite{Champagnat2011,Champagnat2010,Jabin2011,Perthame2007}. 
Those systems are characterised by the existence of at least as many equilibrium states that the number of competiting species (or genotypes) involved. In addition, the dynamics exhibit a competitive exclusion principle which state that the fittest species initially present  will overcome all the other ones.

When the mutations's matrix $(\mu_{ij})$ is  non trivial,  the system \eqref{mutselecedo-eq-sysgen} does not fall into  Hirsch's definition of competitive system and  less results are known.  If  for reasonably smooth interaction functions, the existence of solution of \eqref{mutselecedo-eq-sysgen} defined for all times is not an issue, the existence of a non trivial stationary solution 
and the analysis of the asymptotic behaviour are challenging questions.

Most of the known results concerns either particular interaction functions $\Psi_i$ for which the existence of steady states and their local stability are investigated \cite{Bates2011,Crow1970,Hadeler1981,Hofbauer1985} or for some  ODE's systems \eqref{mutselecedo-eq-sysgen} where   the mutation matrix $\mu_{ij}$ is considered as a small parameter. In the latter the system \eqref{mutselecedo-eq-sysgen} is  seen as a perturbation  of  \eqref{mutselecedo-eq-sysgen2} and analysed using perturbative techniques \cite{Burger1994, Calsina2005,Calsina2007}. 
Recently, there has been an intense activity on continuous trait version of \eqref{mutselecedo-eq-sysgen} where some of the techniques can be used to obtain  the existence of locally stable steady states for  \eqref{mutselecedo-eq-sysgen}  see for example \cite{Barles2008,Calsina2005,Calsina2007,Calsina2012,Canizo,Carrillo2007,Champagnat2011,Champagnat2010,
Desvillettes2008,Diekmann2005,Jabin2011,Perthame2007, Raoul2011,Raoul2012}.
However,  to our knowledge there is no results on global stability of the steady states for systems like\eqref{mutselecedo-eq-sysgen} neither for its continuous trait version.


 

In this direction, our first results concern the systems \eqref{mutselecedo-eq-sysgen} where the competition terms $\Psi_i$ are independent of $i$.  A typical case is $\Psi_i(v)=\sum_{j=1}^Nv_j$  which corresponds to a situation where inter- and intra-species (genotypes) competitions are equals (i.e. blind and uniform competition between variants,  see Lafforgue et al. \cite{Lafforgue2011}).
 In this situation one can show that there exists a unique   positive stationary solution  of \eqref{mutselecedo-eq-sysgen} $\bar v\in \R^N$, which attracts all the trajectories initiated  from any nonnegative and non zero initial data.  
Namely we show that 
\begin{theorem} \label{mutselecedo-thm1}
Assume that the interaction $\Psi_i$ is independent of $i$ and is satisfying  the assumptions  \eqref{mutselecedo-hyp1}--\eqref{mutselecedo-hyp2}, then there exists a unique positive stationary solution $\bar v$ to  \eqref{mutselecedo-eq-sysgen}. Moreover for any nonnegative  initial datum $v(0)$ not identically  zero, the  corresponding solution $v(t)$ of \eqref{mutselecedo-eq-sysgen}  converges to $\bar v$.  
\end{theorem}

A case of particular interest  is when the  interactions $\Psi_{i}$ take the following form $\Psi_i(v)=\sum_{j=1}^Nr_jv_j$.  This particular structure of interaction  was initially introduced on a theoretical ground by Sole et al. \cite{Sole1999} to model the competition between viruses. Recently, this form of interactions has been used to explain experimental results of virus evolution \cite{Cuevas2003,Fabre2012,Miralles2001}.
Sole et al. showed that Eigen's model of molecular quasi-species \cite{Eigen1988} was to a large extent equivalent to the Lotka-Volterra competition equations under this assumption.

For this type of interaction,  we can rewrite the system \eqref{mutselecedo-eq-sysgen}  as follows
 \begin{equation}
 \frac{dv_{i}(t)}{dt}=v_i\left[r_i-\frac{1}{K}\sum_{j=1}^{N}r_{j}v_j\right]+\sum_{j=1}^{N} \mu_{ij}(v_j-v_i).\label{mutselecedo-eq-syspart}
 \end{equation}
For \eqref{mutselecedo-eq-syspart} besides the asymptotic behaviour of the solution obtained as an application of  Theorem \ref{mutselecedo-thm1} we can precise  the speed of convergence to the equilibrium. Furthermore we can  give an estimate of the time to reach near the equilibrium.  
In epidemiology, this type of informations is of practical use  for building tractable nested models. Nested models are a class of model which explicitly links the relationships between processes at different levels of biological organization.They are often used to study the pathogen evolution by linking the disease dynamic of within- and between-host, see Mideo et al. \cite{Mideo2008}. Their formalisation becomes more simple when the within-host pathogen dynamic is faster than the between-hosts dynamic. Indeed, in such cases, using for example slow-fast reduction techniques commonly used in ecology  \cite{Auger2008}, the within-host dynamic can be approximated by its equilibrium, see for example \cite{Fabre2012a}.

More precisely, we show that 

\begin{theorem}\label{mutselecedo-thm2}
Assume that the interactions $\Psi_i$ take the form $\Psi_i(v)=\sum_{j=1}^Nr_jv_j$ then for any nonnegative  initial datum $v(0)$ not identically  zero, the  solution $v(t)$ of \eqref{mutselecedo-eq-syspart} converges exponentially fast to its unique equilibrium. That is to say there exists two positive constants $C_1$ and $C_2$ so that 
$$\|v-\bar v\|_{\infty}\le C_1e^{-C_2t}.$$ 
\end{theorem}

Next, we analyse the situation where the functions $\Psi_i$ are not reduced to a single function. From a biological point of view, this situation could appears as, often, genotypes (species) exhibit particular association patterns (see for example \cite{Takahashi2007} in the case of plant viruses). A way to model such phenomena is to take $\Psi_i$ of the form
$\Psi_i(v)=\sum_{j=1}^N \alpha_{i,j}r_jv_j$
where $\alpha_{i,j}$ are crowing index \cite{Lloyd1967,Hartley2002}. This index is equal to 1 when the 2 species are distributed independently and ranges from 0 (complete avoidance) to a large constant (near overlap) according to patterns of species association.

 In this general context, our first result concerns the existence of a positive stationary solution  for the  system \eqref{mutselecedo-eq-sysgen} assuming we have  the following extra condition 
 
 \begin{align}
 \label{mutselecedo-hyp3} \forall\; i\quad \sum_{j=1}^{N}\mu_{ij}\le \frac{r_i}{2}.
 \end{align}
This condition makes senses in our application framework  since mutation rates (expressed as substitutions per nucleotide per cell infection) range from $10^{-8}$ to $10^{-6}$ for DNA viruses and from $10^{-6}$ to $10^{-3}$ for RNA viruses \cite{Sanjuan2010} whereas growth rates $r_i$ are  of the order of the unit or higher. For example, the overall growth rate of the RNA virus VSV was estimated to $0.6$ virus/hour \cite{Cuevas2005}. Each cell infected  by a single VSV particle produces from 50 to 8000 progeny virus particle/ cell infection \cite{Zhu2009}. 

Under the extra assumption \eqref{mutselecedo-hyp3} we prove  
\begin{theorem}\label{mutselecedo-thm3}
Assume that  $r_i,\Psi_i$ and $\mu_{ij}$ satisfy the assumptions \eqref{mutselecedo-hyp1}--\eqref{mutselecedo-hyp2}. Assume further that \eqref{mutselecedo-hyp3} holds .  Then there exists a positive  stationary solution $\bar v$ to the system \eqref{mutselecedo-eq-sysgen}.
\end{theorem}

Lastly, we obtain the convergence to the equilibria for some general interactions $\Psi_i$.  Namely, we show that
\begin{theorem}\label{mutselecedo-thm4}
Assume that  $r_i,\Psi_i$ and $\mu_{ij}$ satisfy \eqref{mutselecedo-hyp1}--\eqref{mutselecedo-hyp2}. Assume further that \eqref{mutselecedo-hyp3} holds and $\Psi_i(v)=\alpha(v)+\eps\psi_i(v)$ with $\psi_i\in C^1_{loc}$ uniformly bounded and $\alpha \in C^{1}_{loc}$ satisfying \eqref{mutselecedo-hyp1}--\eqref{mutselecedo-hyp2} is so that  $\forall\, x\in \R^N,\; \nabla\alpha (x)>0$.  Then there exists $\eps_0$ so that, for all $\eps \le \eps_0$, the positive  stationary solution $\bar v_\eps$ of \eqref{mutselecedo-eq-sysgen} attracts all the possible trajectories initiated from any non zero and  nonnegative initial data.
\end{theorem}

 Note that $C^1_{loc}$ perturbation of the particular interaction function $\alpha(v)=\sum_{j=1}^Nr_jv_j$ satisfied the assumptions of the above Theorem.



\medskip

\subsection{General remarks} ~

Before going to the proofs,  we want to make some general remarks and comments.

\subsubsection*{Mathematical comments :}
First  as remarked in \cite{Champagnat2010},  we can interpret  the system of  equations \eqref{mutselecedo-eq-syspart}  as a discrete version of the continuous model 
\begin{equation}
\frac{\partial v(t,x)}{\partial t}=v(t,x)\left[r(x)-\frac{1}{K}\int_{\R^n}r(y)v(t,dy)\right]+\int_{\R^n}\mu(x,y)v(t,dy)  -\mu(x)v(t,x) \label{mutselecedo-eq}
\end{equation}
by posing for each sub population corresponding to a typical trait  $x_i \in \R^d$, 
\begin{align*}
&v(t,x)=\sum_{i}^{N}v_i(t)\delta_{x_i}, \; r(x_i)=r_i,\; \mu(x_i,x_j)=\mu_{ij} \text{ and }\mu(x)=\int_{\R^n}\mu(x,y)dy.
\end{align*}

Recently there have been a lot of works dealing with \eqref{mutselecedo-eq} and generalisation of it, see for example \cite{Barles2008,Burger2000,Calsina2005,Calsina2007,Calsina2012,Canizo,Carrillo2007,Champagnat2010,Champagnat2011,
Desvillettes2008,Diekmann2005,Jabin2011, Perthame2007,Raoul2011,Raoul2012} and references therein. 
A large part of the analysis are concerning \eqref{mutselecedo-eq} in absence of mutation  (ie $\mu \equiv 0$) or in the limit $\mu \to 0$.  In the latter case, much  have been done  in developing a constrained Hamilton-Jacobi approach to analyse the long time behaviour of positive solution of this type of models see for instance \cite{Barles2008,Diekmann2005}.  Other approaches based on semigroup theory  have also been developed  to analyse the asymptotic behaviour and local stability of the stationary solution of \eqref{mutselecedo-eq}  see  \cite{Calsina2005,Calsina2007}. Although some of the techniques  developed in this two frameworks  may be adapted to analyse the system \eqref{mutselecedo-eq-sysgen} most of them fail when we try to prove the global stability of the stationary states. 
To tackle this difficulty we construct a set of Lyapunov functionals for the solution of \eqref{mutselecedo-eq-sysgen} when the interactions $\Psi_i$ are independent of $i$.  Properly used these Lyapunov functionals enable us to analyse the asymptotic behaviour of the solution of \eqref{mutselecedo-eq-sysgen} in this particular  situation. From this analysis, we derive some consideration on the  asymptotic behaviour of the solution of \eqref{mutselecedo-eq-sysgen} in a general situation. 
The  Lyapunov functionals are constructed in the spirit of the relative entropy introduced for linear parabolic operators in \cite{Michel2005}.  Even if  our problem is nonlinear such type of relative entropy can still  be constructed. These is an interesting new feature  of nonlinear dissipative system.

It is worth noticing that a similar construction can be made for \eqref{mutselecedo-eq} giving us access to a simple way of analysing the asymptotic behaviour and the global stability of steady solution of  \eqref{mutselecedo-eq}, see \cite{Coville2012a}.

\medskip
Along  some of the   proofs  we notice that  the existence of a steady state in Theorem \ref{mutselecedo-thm1} and  Theorem \ref{mutselecedo-thm3} can be generalised to  situation where the mutation matrix $(\mu_{ij})$ is not  symmetric. 
Indeed, when the interactions  $\Psi_i$ are independent of $i$ (Theorem \ref{mutselecedo-thm1}), the construction of   a unique stationary solution relies only on  the Perron-Frobenius Theorem  which holds true for non symmetric matrices. However, for the general case (Theorem \ref{mutselecedo-thm3}) to obtain the existence of a steady state,    we do require that the condition \eqref{mutselecedo-hyp3} is replaced by
\begin{align}
 \label{mutselecedo-hyp4} \forall\; i\quad \sum_{j=1}^{N}\left(\frac{\mu_{ij}+\mu_{ji}}{2}\right)\le \frac{r_i}{3}.
 \end{align}
 \subsubsection*{Biological comments :}
 First, we emphasize the biological interpretation of our result and particularly
the role of the mutation term. Indeed,   under  biological compatible assumptions concerning the competition $\psi_i$, the mutation matrix $(\mu_{ij})$ and the growth rate $r_i$ we have shown that, in sharp contrast with the classical
results known for the Lotka-Volterra system without mutation, the mutation term deeply changes the dynamics of the system. 
On one hand, the mutation term stabilizes the dynamic
of system by reducing the number of equilibrium up to a single equilibrium and  on
the other hand, the mutation term precludes the competitive exclusion principle to occur. 


Second, we emphasize the biological relevance of relaxing some hypothesis regarding (i) the
monotony of $\Psi_i$ and (ii) the symmetry of the mutation matrix $\mu_{ij}$ to study viral demo-genetics
dynamics in structured hosts. 
To illustrate this point,  let consider 2 patches $p_1$ and $p_2$, and 2 virus genotypes
$v_1$ and $v_2$. The demo-genetics dynamics of the viral population in this system can be modelled by

$$\frac{dw}{dt}=Mw+Rw-\Psi(w)w$$
where $w$ is the vector $(v_{1,p_1},v_{2,p_1},v_{1,p_2},v_{2,p_2})$ and $R,M$ and $\Psi(w)$ are the following matrices 
$R:= (r_i \delta_{ij})$, 

 $$ M:= \begin{pmatrix}  \mu&1-\mu &0 &0 \\
                                              1-\mu &  \mu  & 0  & 0   \\
   0 &0  &\mu& 1-\mu\\
   0&0 &1-\mu&\mu
\end{pmatrix} +\begin{pmatrix}   
-d_1 &0 &d_{2}& 0 \\
0 &-d_1 & 0 & d_{2}\\
 d_{1}&0 &-d_{2}  &0  \\
 0&d_{1}&0 &-d_{2} 
\end{pmatrix} , $$
$$
 \Psi(w):= \begin{pmatrix} r_1v_{1,p_1}+r_2v_{2,p_1} &0&0 &0 \\
    0& r_1v_{1,p_1}+r_2v_{2,p_1} &0 &0 \\
    0&0&r_1v_{1,p_2}+r_2v_{2,p_2}&0 \\
    0&0&0& r_1v_{1,p_2}+r_2v_{2,p_2}  \\
\end{pmatrix}. $$ 


 
In structured hosts some tissue types often act as virus "sources" and others are "sinks", creating an asymmetry in the exchange between patches. In the above example this implies that $d_1 > d_2$ if the patch $p_1$ is a "source" and the patch $p_2$ is a "sink", making the mutation matrix $M$ non-symmetric. 
Moreover, in this example since the competition takes place only inside a given patch, we can check that the two monotone interaction functionals involved  $\Psi_1(v):=r_1v_{1,p_1}+r_2v_{2,p_1}$ and $\Psi_2(v):=r_1v_{1,p_2}+r_2v_{2,p_2}$   do not satisfy the monotone properties \eqref{mutselecedo-hyp1}.
 
 Such type of structure,   inducing asymmetries  in the exchanges and some weak interaction functionals,   are expected to play a role in constraining or facilitating adaptive evolution of viruses in heterogeneous host environment \cite{Orive2005}. It is then relevant to study extension of our results in a context of more general assumptions on $\mu_{ij}$ and $\Psi_i$.

\medskip
\subsubsection*{Organization of the Paper :} In  Section \ref{mutselecedo-section-1}   we start by establishing some preliminaries results  about the system \eqref{mutselecedo-eq-sysgen}.  We derive the relative entropy identities in \ref{mutselecedo-section-re}
Then in Section \ref{mutselecedo-section-specialcase} we analyse in details the system \eqref{mutselecedo-eq-sysgen} for a particular type of interactions and we prove the Theorem \ref{mutselecedo-thm1}.  We show Theorem \ref{mutselecedo-thm2} in Section \ref{mutselecedo-section-caseofinterest}. Finally, we are proving the Theorems \ref{mutselecedo-thm3} and \ref{mutselecedo-thm4} in Sections  \ref{mutselecedo-section-generalcase} and  \ref{mutselecedo-section-genasb}.

\section{ Global facts on the system \eqref{mutselecedo-eq-sysgen}}\label{mutselecedo-section-1}
In this section we establish some useful properties of solution of \eqref{mutselecedo-eq-sysgen}  and  prove the existence of a positive global in time solution of the system \eqref{mutselecedo-eq-sysgen}, that for convenience we rewrite 

 \begin{equation}
 \frac{dv}{dt}=A(\Psi(t))v \label{mutselecedo-eq-sysgenred}
 \end{equation}
 where $A(\Psi(t))$ is the following matrix:
 $$ A:= \begin{pmatrix} (r_1-\frac{1}{K}\Psi_1(t)) -\mu_1+\mu_{11}  & &\mu_{ij} \\
     & \ddots &  \\
   \mu_{ij} &  & (r_{N}-\frac{1}{K}\Psi_N(t))-\mu_{N}+ \mu_{NN}\\
\end{pmatrix}, $$ 
with $\Psi_i(t):=\Psi_i(v(t))$ and $\mu_{i}=\sum_{j=1}^{N}\mu_{ij}$.
Since the function $\Psi_i$ are locally Lipschitz  the local existence of a solution of \eqref{mutselecedo-eq-sysgenred} is a straightforward application of the Cauchy-Lipschitz Theorem. To obtain a global solution starting with non-negative  initial data, we need more \textit{ a priori } estimates on the solutions $(v_i)_{i\in \{1,\ldots,N\}}$. Let  us  first show that the system \eqref{mutselecedo-eq-sysgenred} preserve the positivity.

\subsection{Positivity}
\begin{lemma}[Positivity]\label{mutselecedo-lem1}
The system \eqref{mutselecedo-eq-sysgenred} preserves the positivity. That is to says that if $v(0)$ is a nonnegative initial value, then the solution $v$ of \eqref{mutselecedo-eq-sysgenred} stays  nonnegative. Moreover,  for all $i$ $v_i(t)>0$ for all times $t\in (0,T]$.
\end{lemma}
\dem{Proof:}
We argue by contradiction and assume that $v$ changes its sign.  Let $t_1\in (0,T]$ be the first time where   $v_i(t_1)=0$ for some $i$, $v(t)>0$ for all $t\in(0,t_1)$ and there exists $t>t_1$ so that $v(t)<0$. $t_1>0$ is well defined since $v$ changes its sign and  $v_i(0)\ge 0$ for all $i$, $v_i(0)>0$ for some $i$ and $\lim_{t \to 0^+}\frac{dv_i(t)}{dt}>0$.
 Now  at the time $t_1$, we have
$$0\ge \frac{dv_i(t_1)}{dt}=\sum_{j}\mu_{ij}v_j(t_1)\ge 0.$$ 
Therefore for all $i$  $v_i(t_1)=0=\frac{dv_i(t_1)}{dt}$ and  by the Cauchy-Lipschitz Theorem we deduce that for all $i$ and all 
$t\in [t_1,T]$ $v_i(t)=0$. Thus $v\ge 0$ in $(0,T]$  which contradicts that $v$ changes its sign.

Next we show  that $v_i(t)>0$ for all times $t\in (0,T]$. From the above argumentation, $v$ is nonnegative for all times. To show that $v$ stays positive for all time,  we can see that it is enough to show that $\n(t)=\sum_{i=1}^{N}v_j(t)>0$ for all times $t\in [0,T]$. Indeed, if there is  $t_1\in (0,T]$  so that $t_1$ is the first time where for some $i$,  $v_i(t_1)=0$ then arguing as above we see that for all $i$ and $\t\ge t_1$ $v_i(t)=0$ so $\n(t)=0$ for all $t\ge t_1$.

Now,  let us denote $Q(0,1)$ the following set $$Q(0,1):=\left\{x\in\R^N\;| \sum_{i=1}^N|x_i|\le 1\right\}.$$
  Since $\Psi_i$ is locally Lipschitz, we have on $Q(0,1)$,  
 \begin{equation}\label{mutselecedo-lip1}\forall \, i, \exists \kappa_i \; \text{ so that }\; \Psi_i(x)\le \kappa_i\sum_{i=1}^N|x_i|.
 \end{equation}

 Let $\kappa_0:= \sum_{i=1}^{N}\kappa_i$ and $\delta$ a positive constant  so that $\delta<\min\left\{1,\frac{\n(0)}{2}, \frac{r_{min}}{2\kappa_0}\right\}$. We will prove that $\n\ge \delta$ for all times in $(0,T]$.  We argue by contradiction and assume there exists $t_1<t_2$ so that $\n(t_1)=\delta$ and $\n(t)\le \delta$ in $(t_1,t_2]$. 
 Using  that the $v_i$ are non negative functions and \eqref{mutselecedo-lip1}, from our assumption  we deduce that on  $[t_1,t_2]$  and  for all $i$,  
$$\frac{dv_i}{dt}\ge v_i\left(r_{min}-\kappa_i\n\right)+\sum_{j=1}^{N}\mu_{ij}(v_j-v_i).$$ 
By summing over all the possible $i$  we end up with  the following equation
\begin{align}
&\frac{d\n}{dt}\ge \n(r_{min}-\n\kappa_0),\label{mutselecedo-eq-lpos1}\\
&\n(t_1)=\delta.
\end{align}
 Using the  comparison principle on the latter equation, we achieve on $[t_1,t_2],$ $\n\ge \tilde \n$ where $\tilde \n$ is the solution of the logistic type equation \eqref{mutselecedo-eq-lpos1}. By construction    we have  $\tilde \n>\delta$ for all times $t>t_1$. Thus  we get the contradiction $\delta<\n\le \delta$. Hence   $\n>\delta$ for all $t \in (0,T]$.
 
\fdem
\begin{remark}\label{mutselecedo-rem-lem1}
Of the above proof, we  also have  a bound from below for $\n$. Namely, we have for all $t$
$$\n(t)\ge \min\left\{1,\frac{\n(0)}{2}, \frac{r_{min}}{2\kappa_0}\right\}.$$ 
\end{remark}
\subsection{Existence of a global solution}~\\

Next we show the following estimate an

\begin{lemma}\label{mutselecedo-lem2}
Let $v(t)$ be a solution of the Cauchy problem \eqref{mutselecedo-eq-sysgenred}, then there exists two constants $C_0$ and $C_1$ and a positive vector $v_p$ so that $v(t)\le C_0e^{C_1t}v_p$
\end{lemma} 

\dem{Proof:}
From Lemma \ref{mutselecedo-lem1} we know that the $v_i$ are non negative functions. Therefore by \eqref{mutselecedo-eq-sysgenred} we can see that the $v_i$ satisfy
\begin{equation}
 \frac{dv_i}{dt}\le r_iv_i +\sum_{j=1}^{N}\mu_{ij}(v_j-v_i) \label{mutselecedo-eq-lem3}.
 \end{equation}
  Let us denote $M$ and $R$ the two following matrices   
  $$ M:= \begin{pmatrix}  -\mu_1+\mu_{11}  & &\mu_{ij} \\
     & \ddots &  \\
   \mu_{ij} &  & -\mu_{N}+ \mu_{NN}\\
\end{pmatrix}, \qquad R:= \begin{pmatrix} r_1  & &0 \\
     & \ddots &  \\
   0 &  & r_{N}\\
\end{pmatrix}. $$ 

  We can then rewrite the inequalities \eqref{mutselecedo-eq-lem3} as follows
 \begin{equation}
 \frac{d v}{dt}\le (R+M)v.\label{mutselecedo-eqlem4}
 \end{equation}
 By choosing $\bar \mu=\sup_{i\in\{1,\ldots,N\}}\mu_i$, we can see that  $R+M+\bar \mu$ is nonnegative matrix. Since $R+M+\bar \mu Id$ is also irreducible, by the Perron-Frobenius Theorem $R+M+\bar \mu Id$ posses a unique principal eigenpair $(\nu_p,v_p)$  so that  $v_p$ is  a  positive vector, i.e. there exists $(\nu_p, v_p)$ so that $v_p> 0$ and
 $$
 (R+M+\bar \mu Id)v_p=\nu_p v_p.
 $$
Note that $e^{(\nu_p-\bar \mu)t}\|v(0)\|_{\infty}v_p$ satisfies 
$$ \frac{d u}{dt}= (R+M)u.$$
So, by the Cauchy-Lipschitz Theorem, from \eqref{mutselecedo-eqlem4} we deduce   that $v\le e^{(\nu_p-\bar \mu)t}\|v(0)\|_{\infty}v_p$ for all times $t$.
\fdem

The existence of a global in time solution for the system \eqref{mutselecedo-eq-sysgenred} is then a consequence of the  Cauchy-Lipschitz Theorem and the above a priori estimates.

\begin{remark}
By adapting some ideas  in \cite{Perthame2007} and for a particular type of $\Psi_i$,   we can derive an explicit formula for the solution  of \eqref{mutselecedo-eq-sysgenred}. Indeed when  the interaction terms take the form   $\Psi_{i}(v)=\alpha(v)=\sum_{j=1}^{N}\alpha_j v_j$ for all $i$ then the solution $v_i(t)$ can be expressed by the formula
\begin{equation}
v_i(t)=\frac{\left(e^{(R+M)t}v(0)\right)_i}{1+\sum_{j=1}^N\alpha_j \int_{0}^t\left(e^{(R+M)s}v(0)\right)_j\, ds}.
 \end{equation}

 
To obtain this formula we argue as in \cite{Perthame2007} (Chapter 2 Section 2.1 ) and  we start by introducing the functions 
$$\alpha(t):=\sum_{j=1}^N\alpha_jv_j\quad \text{ and }\quad V_i(t):=e^{\int_{0}^t\alpha(s)\,ds}v_i(t).$$
We remark that the $V_i$ satisfy the linear equation
  \begin{equation}
 \frac{dV_{i}(t)}{dt}=r_iV_i + \sum_{j=1}^{N} \mu_{ij}(V_j-V_i) \label{mutselecedo-eq-solexp1}.
 \end{equation}
Thus $V_i(t):=\left(e^{(R+M)t}v(0)\right)_i$ and $v_i(t)$ is implicitly given by the formula
$$ v_i(t)= e^{-\int_{0}^t\alpha(s)\,ds}\left(e^{(R+M)t}v(0)\right)_i.$$
Now let us evaluate the term $ e^{-\int_{0}^t\alpha(s)\,ds}$. By differentiating  $e^{-\int_{0}^t\alpha(s)\,ds}$ one has
$$
\frac{d}{dt}(e^{\int_{0}^t\alpha(s)\,ds}) =  \alpha(t)e^{\int_{0}^t\alpha(s)\,ds}
=\sum_{j=1}^N\alpha_j V_j(t).
$$
Therefore we have
$$
e^{\int_{0}^t\alpha(s)\, ds} =1+\int_{0}^t\sum_{j=1}^N\alpha_j V_j(s)\,ds 
=1+\sum_{j=1}^N\alpha_j \int_{0}^t\left(e^{(R+M)s}v(0)\right)_j\, ds.
$$

Hence
$$
v_i(t)=\frac{\left(e^{(R+M)t}v(0)\right)_i}{1+\sum_{j=1}^N\alpha_j \int_{0}^t\left(e^{(R+M)s}v(0)\right)_j\, ds}.
$$
\fdem
\end{remark}
\medskip

\section{Relative Entropy Identities.}\label{mutselecedo-section-re}~\\

Here, we prove the following general principle which give us access to  some useful identifies that we constantly use along this paper.
\begin{theorem}\label{mutselecedo-thm-genid}
Assume that $r_i,(\mu_{ij}), \Psi_i$ satisfies \eqref{mutselecedo-hyp1}-- \eqref{mutselecedo-hyp2}.  Let $v$  and $\bar v$ be respectively a positive  solution  and a positive stationary solution   of \eqref{mutselecedo-eq-sysgen}, and let $H$ be a smooth (at least $C^1$)  function.
Then the function $\h(v):=\sum_{i=1}^N\bar v_i^2H\left(\frac{v_i}{\bar v_i }\right)$ satisfies 
\begin{equation}
\frac{d \h(t)}{dt}=-\D(v)+\frac{1}{K}\sum_{i=1}^N \bar v_iH^{\prime}\left(\frac{v_i}{\bar v_i}\right) \Gamma_iv_i
\end{equation} 
where 
\begin{align*}
 &\D(v):=\sum_{i,j=1}^N\mu_{ij}\bar v_i \bar v_j \left[H\left(\frac{v_j}{\bar v_j }\right)-H\left(\frac{v_i}{\bar v_i }\right)\right]
+\sum_{i,j=1}^N\mu_{ij}\bar v_i \bar v_j H^{\prime}\left(\frac{v_i}{\bar v_i }\right)\left[\frac{v_i}{\bar v_i }-\frac{v_j}{\bar v_j }\right]\\
&\Gamma_i:=\Psi_i(\bar v) -\Psi_i(v)
\end{align*}
\end{theorem}

\dem{Proof:}

By \eqref{mutselecedo-eq-sysgen}, for all $i$ we have 

 \begin{equation}
 \frac{d v_i}{dt}= \left(r_i v_i -\frac{1}{K}\Psi_i(\bar v)v_i  +\sum_{j=1}^{N}\mu_{ij}(v_j-v_i)\right) +\frac{1}{K}\Gamma_i v_i
 \end{equation}
 
Using that $\bar v$ is a stationary solution, we have for all $i$
$$(r_i-\frac{1}{K}\Psi_i(\bar v) )\bar v_i=- \sum_{j=1}^{N}\mu_{ij}(\bar v_j-\bar v_i),$$
and we can rewrite the above equation  as follows
  $$
  \frac{d v_i}{dt}=\sum_{j=1}^{N}\mu_{ij}\left(v_j-\frac{\bar v_jv_i}{\bar v_i}\right)+\frac{1}{K}\Gamma_iv_i.
  $$
  By multiplying the above equality by $\bar v_iH^{\prime}\left(\frac{v_i}{\bar v_i}\right)$ and by summing over all $i$ we achieve 
$$\sum_{i=1}^N \bar v_iH^{\prime}\left(\frac{v_i}{\bar v_i}\right) \frac{d v_i}{dt}=\frac{1}{K}\sum_{i=1}^N \bar v_iH^{\prime}\left(\frac{v_i}{\bar v_i}\right) \Gamma_iv_i +\sum_{i=1}^N \bar v_iH^{\prime}\left(\frac{v_i}{\bar v_i}\right)\sum_{j=1}^{N}\mu_{ij}\left(v_j-\frac{\bar v_jv_i}{\bar v_i}\right). $$
Thus we have
$$\frac{d \h(t)}{dt}=\frac{1}{K}\sum_{i=1}^N \bar v_iH^{\prime}\left(\frac{v_i}{\bar v_i}\right) \Gamma_iv_i -\sum_{i,j=1}^N \mu_{ij}\bar v_i \bar v_jH^{\prime}\left(\frac{v_i}{\bar v_i}\right)\left(\frac{v_i}{\bar v_i}-\frac{v_j}{\bar v_j}\right). $$
Hence we have 
$$\frac{d \h(t)}{dt}=\frac{1}{K}\sum_{i=1}^N \bar v_iH^{\prime}\left(\frac{v_i}{\bar v_i}\right) \Gamma_iv_i-\D(v). $$
since by symmetry of $\mu_{ij}$, $$\sum_{i,j=1}^N\mu_{ij}\bar v_i \bar v_j \left[H\left(\frac{v_j}{\bar v_j }\right)-H\left(\frac{v_i}{\bar v_i }\right)\right]=0.$$
\fdem

As immediate corollary  of the Theorem \ref{mutselecedo-thm-genid},  we have the following identities that we constantly used along this paper,
\begin{lemma}\label{mutselecedo-lem-general-id}
Assume that $r_i,(\mu_{ij}), \Psi_i$ satisfies \eqref{mutselecedo-hyp1}-- \eqref{mutselecedo-hyp2}.  Let $v$  and $\bar v$ be respectively a positive  solution  and a positive stationary solution   of \eqref{mutselecedo-eq-sysgen}, then we have 
\begin{itemize}
 \item[(i)] $$\frac{d}{dt}\left(\sum_{i=1}^{N}v_i\bar v_i \right)=\frac{1}{K}\sum_{i=1}^N(\Psi_i(\bar v)-\Psi_i(v))v_i\bar v_i$$
 \item[(ii)] $$\ \frac{d}{dt}\left(\sum_{i =1}^{N}v_i^2\right)=-\sum_{i,j =1}^{N}\mu_{ij} \bar v_i \bar v_j\left( \frac{v_j}{\bar v_j}-\frac{v_j}{\bar v_j}\right)^2+\frac{2}{K}\sum_{i =1}^{N}v_i^2\left(\Psi_i(\bar v)-\Psi_i(v)\right).$$
\end{itemize} 
\end{lemma}
\dem{Proof:}
$(i)$ and $(ii)$ can be obtained straightforwardly  from  the Theorem \ref{mutselecedo-thm-genid}. Indeed,  
by using $H(s)=s$ in Theorem \ref{mutselecedo-thm-genid}, and by observing that 
from the symmetry of $\mu_{ij}$,  
$$\sum_{i,j=1}^N\mu_{ij}\bar v_i \bar v_j \left[\frac{v_i}{\bar v_i }-\frac{v_j}{\bar v_j }\right]=0,$$
we easily get
$$\frac{d \h(t)}{dt}=\frac{1}{K}\sum_{i=1}^N \Gamma_i \bar v_iv_i, $$
which proves $(i)$ since for the function $H(s)=s$, we have $\h(t)=\sum_{i=1}^N\bar v_i v_i$.
We can obtain $(ii)$ in a similar way   by using the function $H(s)=s^2$ in the Theorem \ref{mutselecedo-thm-genid} and by observing that 
 $$\sum_{i,j=1}^N\mu_{ij}\bar v_i \bar v_j \left(2\frac{v_i}{\bar v_i }\right)\left[\frac{v_i}{\bar v_i }-\frac{v_j}{\bar v_j }\right]=\sum_{i,j=1}^N\mu_{ij}\bar v_i \bar v_j \left[\frac{v_i}{\bar v_i }-\frac{v_j}{\bar v_j }\right]^2. $$


 \section{The special case $\Psi_{i}(.)=\alpha(.)$:}\label{mutselecedo-section-specialcase}
In this section we analyse in details the asymptotic behaviour of a positive solution of \eqref{mutselecedo-eq-sysgen} when the competition  functional $\Psi_i$ is  independent of $i$, i.e for all $i, \Psi_i(v)=\alpha(v)$ where $\alpha$ is a function from $\R^N\to \R$ which satisfies \eqref{mutselecedo-hyp1} --\eqref{mutselecedo-hyp2}. 
As we expressed in Theorem \ref{mutselecedo-thm1} that we recall below, in this situation the system \eqref{mutselecedo-eq-sysgen} has a unique positive stationary solution which attracts all the trajectories initiated from any nonnegative and non zero initial data. More precisely, we prove 
 
 \begin{theorem}
 Assume that for all $i$ the  function   $\Psi_i(.)=\alpha(.)$, then there exists a unique stationary solution $\bar v$ of \eqref{mutselecedo-eq-sysgenred}. Moreover, for all nonnegative and non zero initial datum $v(0)$, the corresponding solution $v(t)$ converges to $\bar v$.  
 \end{theorem} 
 
 To prove the Theorem, we first analyse the existence of stationary solution

\medskip

\subsection{Study of the existence of equilibria}\label{mutselecedo-ss-equilibria}~\\
Recall that we look for a stationary solution of   
  
  \begin{equation}
 \frac{dv}{dt}=A(\alpha(t))v. \label{mutselecedo-eq-sysred}
 \end{equation}

Therefore if there exists a stationary equilibria for the system \eqref{mutselecedo-eq-sysred} the $v_i$ must satisfies the following equations: 
\begin{align}
&A(\bar \alpha)v=0\label{mutselecedo-eq-equilibria1}\\
&\bar \alpha=\alpha(v)\label{mutselecedo-eq-equilibria2}
\end{align}
 where $A(\bar \alpha)$ is the matrix 
 $$ A(\bar \alpha):= \begin{pmatrix} (r_1-\frac{\bar \alpha}{K}) -\mu_1+\mu_{11}  & &\mu_{ij} \\
     & \ddots &  \\
   \mu_{ij} &  & (r_{N}-\frac{\bar \alpha}{K})-\mu_{N}+ \mu_{NN}\\
\end{pmatrix}. $$ 
 
Note that we can rewrite the matrix $A(\bar \alpha):=\left(R-\left(\frac{\bar \alpha}{K}\right)Id+M\right)$ where $R$ and $M$ are  matrices defined in Section \ref{mutselecedo-section-1}.

  Therefore, a solution $v$ to \eqref{mutselecedo-eq-equilibria1} is a solution to 
  \begin{equation}
  (M+R)v=\left(\frac{\bar\alpha}{K}\right) v.\label{mutselecedo-eq-equilibria3}
  \end{equation}
  Let us now establish some important property of the equilibrium.
  \begin{lemma}
  If $v$ is a nonnegative stationary solution of \eqref{mutselecedo-eq-sysred}, then either $v\equiv 0$ or $v>0$ (i.e $\forall \, i , \, v_i>0$).
  \end{lemma}
  \dem{Proof:}
  First, we observe that $0$ is a solution of the problem \eqref{mutselecedo-eq-equilibria1}.  Now let $v$ be nonnegative stationary  solution  of \eqref{mutselecedo-eq-sysred}  so that $v_i=0$ for some $i$. From \eqref{mutselecedo-eq-equilibria3} we see that 
  $$(M+R v)_i=\left(\frac{\bar \alpha}{K}\right)v_i= 0.$$
  Therefore $$0=r_iv_i+\sum_{j=1}^{N}\mu_{ij}(v_j-v_i)=\sum_{j=1}^{N}\mu_{ij}v_j.$$ 
  Thus $v_j=0$ for all $j$ since by assumption $v_j\ge 0$ and $(\mu_{ij})$ is irreducible.
Hence any nonnegative stationary solution is either positive or the zero solution. 
\fdem  
  
Observe that from the above Lemma and  \eqref{mutselecedo-eq-equilibria3}, one can see that $v$ is a positive eigenvector of the matrix $M+R$ associated to the eigenvalue $\frac{\bar\alpha}{K}$. Now we are in position to prove that there exists a unique positive stationary solution to  \eqref{mutselecedo-eq-equilibria3}.
  
 \begin{lemma}\label{mutselecedo-lem-steady}
 There exists a unique $(\bar\alpha,\bar v)$ solution of the equations  \eqref{mutselecedo-eq-equilibria1} and \eqref{mutselecedo-eq-equilibria2}. Moreover $v$ satisfies $\alpha(v)=\nu_p$.
 \end{lemma} 
 \dem{Proof:}
 
By choosing $\bar \mu=\sup_{i\in\{1,\ldots,N\}}\mu_i$, we can see that  $R+M+\bar \mu$ is nonnegative matrix. Since $R+M+\bar \mu Id$ is also irreducible, by the Perron-Frobenius Theorem $R+M+\bar \mu Id$ posses a unique principal eigenpair $(\nu_p,v_p)$  so that  $v_p$ is  a  positive vector, i.e. there exists $(\nu_p, v_p)$ so that $v_p> 0$ and
 \begin{equation}
 (R+M+\bar \mu Id)v_p=\nu_p v_p.\label{vp}
 \end{equation}
 Moreover, the linear subspace associated to the eigenvalue $\nu_p$ is one dimensional \cite{Zeidler1986}. So without any loss of generality  we can assume that  $\sum_{i=1}^N(v_p)_i^2=1$.
  From the equation \eqref{vp} we deduce that the vector $v_p$ is a positive eigenvector of the matrix $M+R$ associated with the eigenvalue $\lambda_p:=(\nu_p-\bar\mu)$. By construction one can see that $\lambda_p$ is the unique eigenvalue of the matrix $M+R$ associated with a positive eigenvector.
  A quick computation shows that $\lambda_p=(\nu_p-\bar \mu)>0$. Indeed, if not we have 
  $$ (R+M)v_p\le 0.$$
  Thus for all $i\in \{1,\ldots,N\}$ we have 
  
  $$r_i(v_p)_i+\sum_{j=1}^{N}\mu_{ij}\big((v_p)_j-(v_p)_i\big)\le 0.$$ 
  Let $(v_p)_{i_0}:=\min_{i\in\{1,\ldots,N\}}(v_p)_i$ then for $(v_p)_{i_0}$ we have 
  $$ \sum_{j=1}^{N}\mu_{ij}\big((v_p)_j-(v_p)_{i_0}\big)\ge 0.$$  Since $R$ is a positive matrix we achieve the contradiction
  $$0<r_{i_{0}}(v_p)_{i_0}+ \sum_{j=1}^{N}\mu_{ij}\big((v_p)_j-(v_p)_{i_0}\big)\le 0.$$  

 Now from \eqref{mutselecedo-eq-equilibria3} we deduce that  there exists an unique positive $\bar \alpha$ so that $\frac{\bar \alpha}{K}=\lambda_p$.
 Let us now construct our solution. Note that for any $\mu \in \R$, the vector $\mu v_p$ is also a solution to \eqref{mutselecedo-eq-equilibria3} with the eigenvalue $\lambda_p$.  So to obtain a solution $\bar v$  to \eqref{mutselecedo-eq-equilibria1} and  \eqref{mutselecedo-eq-equilibria2} we only have to adjust $\mu$ in such a way that $\alpha(\mu v_p) =\bar \alpha$.  This is always possible for a unique $\mu$ since  $\alpha(0)=0$, $\lim_{\mu\to \infty}\alpha(\mu v_p)=+\infty$ and $\alpha$ is an increasing function.

 \fdem

 \subsection{Convergence to the unique equilibrium}\label{mutselecedo-ss-asb}~\\
 Let us look at the convergence of $v(t)$ toward its equilibrium. Let us first establish some useful identities.
 
 \begin{lemma}\label{mutselecedo.lem-ener1}
 Let us denote $(\bar \alpha,\bar v)$ the stationary solution constructed above. Let $v$ be a solution of \eqref{mutselecedo-eq-sysred} then $v$ satisfies the following identities  
\begin{itemize}
\item[(i)]  $$ \frac{d}{dt}\left(\sum_{i =1}^{N} v_i \bar v_i\right)=\frac{1}{K}(\bar \alpha -\alpha(t))\left(\sum_{i =1}^{N}v_i\bar v_i\right).$$
  
  \item[(ii)]$$ \frac{d}{dt}\sum_{i =1}^{N}\left( v_i \right)^2=-\sum_{i,j =1}^{N}\mu_{ij} \bar v_i \bar v_j\left( \frac{v_j}{\bar v_j}-\frac{v_j}{\bar v_j}\right)^2+\frac{2}{K}(\bar \alpha -\alpha(t))\sum_{i =1}^{N}\left(v_i\right)^2.$$
\end{itemize} 
 \end{lemma}
 
 \dem{Proof:}
Since here for all $i, \Psi_i=\alpha$, from Lemma \ref{mutselecedo-lem-general-id} we deduce that 

\begin{align*}
&\frac{d}{dt}\left(\sum_{i =1}^{N} v_i \bar v_i\right)=\frac{1}{K}\left(\sum_{i =1}^{N}(\alpha(\bar v) -\alpha(v(t)))v_i\bar v_i\right),\\
 &\frac{d}{dt}\sum_{i =1}^{N}\left( v_i \right)^2=-\sum_{i,j =1}^{N}\mu_{ij} \bar v_i \bar v_j\left( \frac{v_j}{\bar v_j}-\frac{v_j}{\bar v_j}\right)^2+\frac{2}{K}(\bar \alpha -\alpha(t))\sum_{i =1}^{N}(\alpha(\bar v)-\alpha(v) )\left(v_i\right)^2.
\end{align*} 
 Thus $(i)$ and $(ii)$ hold true.

  \fdem

From the above Lemma we can derive a  useful Lyapunov functional.
\begin{lemma} \label{mutselecedo.lem-liap}
Let $\bar v$ be  the positive stationary solution of \eqref{mutselecedo-eq-sysred}.  For any positive solution $v$ of \eqref{mutselecedo-eq-sysred}, let  us denote  $\e$, $\beta$ and $F$ the following quantities  
 $$\e(v):=\sum_{i=1}^{N}\left(v_i\right)^2,\quad \beta(v):=\sum_{i=1}^{N}v_i\bar v_i \quad F(v):=\log\left[\frac{\e(v)}{(\beta(v))^2}\right].$$
 Then $F\ge \log(\frac{1}{\sup_{i}\bar v_i})$ and for any  positive solution $v$  of \eqref{mutselecedo-eq-sysred} we have 
 $$\frac{d}{dt} F(v)=-\frac{1}{\e(v)}\sum_{i,j=1}^{N}\mu_{ij}\bar v_j\bar v_i\left(\frac{v_j}{\bar v_j}-\frac{v_i}{\bar v_i}\right)^2<0.$$
 \end{lemma}

\dem{Proof:}
First let us show that $F$ is bounded from below. By construction and using a standard convexity argument, we see that  $\beta^2(v)\le (\sup_{i}\bar v_i)\e(v)$. So from the monotonicity of the $\log$, we conclude that $F(v)\ge \log(\frac{1}{\sup_{i}\bar v_i})$.

Now let us show that  $F$ is non increasing. To do so let  us compute $ \frac{d}{dt} F(v)$.
From the definition of $F$   we see that 
$$ \frac{d}{dt} F(v)=\frac{d_t \e(v)}{\e(v)} -2 \frac{d_t \beta(v)}{\beta(v)},$$
where $d_t$ denotes $\frac{d}{dt}$.
By Lemmas \ref{mutselecedo.lem-ener1}  we have 
\begin{align*}
&\frac{d_t \e(v)}{\e(v)}=-\frac{1}{\e(v)}\sum_{i,j=1}^{N}\mu_{ij}\bar v_j\bar v_i\left(\frac{v_j}{\bar v_j}-\frac{v_i}{\bar v_i}\right)^2+\frac{2}{K}(\bar \alpha -\alpha(t))\\
&\frac{d_t \beta(v)}{\beta(v)}=\frac{1}{K}(\bar \alpha -\alpha(t)).
\end{align*}
Thus, 
$$\frac{d}{dt} F(v)=-\frac{1}{\e(v)}\sum_{i,j=1}^{N}\mu_{ij}\bar v_j\bar v_i\left(\frac{v_j}{\bar v_j}-\frac{v_i}{\bar v_i}\right)^2 \le 0.$$

\fdem

Next we derive  some \textit{ a priori } estimates on the solution $v$ of \eqref{mutselecedo-eq-sysred} from the previous Lemmas. Namely, we show that
\begin{lemma} \label{mutselecedo-lem-esti}
 Let $v$  be a positive solution   of \eqref{mutselecedo-eq-sysred}, then  there exists $C_1$ so that  
 $$\e(v)+\beta(v)\le C_1.   $$ 
\end{lemma}
\dem{Proof:}
Observe that from Lemma \ref{mutselecedo.lem-liap}, to obtain the bound it is sufficient to have a uniform bound on $\beta$. Indeed,
since $F$ is decreasing in time, we have for all times 
$$\e(v)\le \beta^2(t)\frac{\e(v(0))}{\beta^2(v(0))}.$$

Now recall that by Lemma \ref{mutselecedo.lem-ener1}, $\beta(v)$ satisfies the following equation
$$\frac{d \beta(t)}{dt}=\frac{1}{K}(\bar \alpha -\alpha(t))\beta(t).$$

Since $\alpha$ satisfies the assumptions \eqref{mutselecedo-hyp1}--\eqref{mutselecedo-hyp2},  there exists $R_\alpha, c_\alpha, k_\alpha$ so that for all $x\in \R^{N,+}\setminus Q_{R_\alpha}(0)$,
\begin{equation}
c_{\alpha}\left(\sum_{i=1}^Nx_i\right)^{k_{\alpha}}  \le \alpha(x). \label{mutselecedo-eq-estialpha} 
\end{equation}

Assume that for some $t>0, v(t) \in  \R^{N,+}\setminus \bar Q_{R_\alpha}(0)$ otherwise we are done since 
$$\beta(v)\le \left(\max_{i\in \{1,\ldots, N\}}\bar v_i\right) \sum_{j=1}^N|v_j|.$$
Let   $\Sigma$ be
$$
\Sigma:=\{t\in \R^+\, |\, \n(t)>R_{\alpha}\}.
$$
 From \eqref{mutselecedo-eq-estialpha} and using that the  $v_i$ are non negative,   we see that  for all  $ t  \in\Sigma$
 $$ \alpha(t)\ge c_\alpha\n(t)^{k_\alpha} \ge c_\alpha\left(\frac{1}{\bar v_{max}}\right)^{k_\alpha}\beta(t)^{k_\alpha} , $$
where $\bar v_{max}:=\max_{i \in \{1,\ldots, N\}}\bar v_i$. Therefore, on $\Sigma$ we have
 $$\frac{d \beta(t)}{dt}\le \frac{1}{K}\left(\bar \alpha -\tilde c_0\beta^{k_{\alpha}}(t)\right)\beta(t).$$
 Using the logistic character of the above equation, we can check that 
 $$\beta(v(t))\le \sup\left\{\beta(v(0)),\max_{x\in \bar Q_{R_{\alpha}}(0)}\beta(x),\left(\frac{\bar \alpha}{\tilde c_0}\right)^{\frac{1}{k_{\alpha}}}\right\}.$$

\fdem

We are now in position to prove the convergence of $v$ toward its equilibrium.

\begin{lemma}\label{mutselecedo-lem-asb}
Let $(\bar v_i)_{i=1\ldots N}$ be the unique stationary solution of \eqref{mutselecedo-eq-sysred}. Then for any non negative initial datum $(v_i(0))_{i=1\ldots N}$ not identically  zero,  the corresponding solution $(v_i(t))_{i=1 \ldots N}$ of \eqref{mutselecedo-eq-sysred} converges to $(\bar v_i)_{i=1\ldots N}$ as $t$ goes to infinity. 
\end{lemma}
\dem{Proof:}
For simplicity we denote $<,>$  the standard scalar product in $\R^N$.

 Now, since $\bar v \neq 0$ and for all times $t$, $v(t)=(v_1,\ldots,v_N)$ is a vector of $\R^N$, we can write $v(t):=\lambda(t)\bar v +h(t)$ with for all $t,$ $<h(t),\bar v> =0$. 
Substituting $v$ by this decomposition in \eqref{mutselecedo-eq-sysred}, it follows that 
\begin{align}
\lambda'(t)\bar v +\frac{dh(t)}{dt}&=\lambda(t)A(\alpha(t))\bar v +A(\alpha(t))h,\\
&=\frac{1}{K}(\bar \alpha -\alpha(t))\lambda(t)\bar v+A(\alpha(t))h. \label{mutselecedo-eq-sysred-decomp1}
\end{align} 
Therefore, we have 

$$<\lambda'(t)\bar v +\frac{dh(t)}{dt},h>= <\frac{1}{K}(\bar \alpha -\alpha(t))\bar v+A(\alpha(t))h,h>.$$
Thus 
$$<\frac{dh(t)}{dt},h>= \frac{1}{2}\frac{d\e(h)}{dt}=<A(\alpha(t))h,h>.$$
By following the computation developed for the proof of (ii) in Lemma \ref{mutselecedo.lem-ener1}, we see that 
$$\frac{d\e(h)}{dt}=  -\sum_{i,j =1}^{N}\mu_{ij} \bar v_i \bar v_j\left( \frac{h_j}{\bar v_j}-\frac{h_j}{\bar v_j}\right)^2+\frac{2}{K}(\bar \alpha -\alpha(t))\e(h).$$

Since $\e(h)\ge 0$ for all times, we will analyse separately  two  situations: Either $\e(h(t))>0$ for all times $t$ or  there exists $t_0\in \R$ so that $\e(h(t_0))=0$. In the latter case,  from the above equation we see that  we must have $\e(h(t))=0$ for all $t\ge t_0$ and so for all $t\ge t_0$,  we  must have $v(t)=\lambda(t)\bar v$. Hence from \eqref{mutselecedo-eq-sysred-decomp1} we are reduced to analyse the following ODE equation
$$\lambda'(t)=\frac{\lambda(t)}{K}(\bar \alpha -\tilde \alpha(\lambda(t))) $$
 where $\tilde \alpha $ is the increasing locally Lipschitz function defined by $\tilde \alpha(s):=\alpha(s\bar v)$. Note that since $\lambda(t)<\bar v,\bar v>=\beta(v)> 0$, we have $\lambda(t)\ge 0$ for all times $t$. 
The above ODE is of logistic type with non negative initial datum therefore by a standard argumentation we see that  $\lambda(t)$ converges to  $\bar \lambda>0$ where $\bar \lambda$ is the unique solution of  $\tilde \alpha (\bar \lambda )=\bar \alpha$. By construction we have $\tilde \alpha(1)=\bar \alpha$, so we deduce that $\bar \lambda=1$.
Hence, in this situation,   $v$ converges  to $\bar v$ as time goes to infinity.
 
In the other situation,  $\e(h(t))>0$ for all $t$ and we claim that 
\begin{claim} \label{mutselecedo-cla-energy}
 $\e(h(t))\to 0$ as $t\to +\infty$.
 \end{claim}
 Assume the Claim holds true then we can conclude  the proof by arguing as follows.
From  the decomposition $v(t)=\lambda(t)\bar v +h(t)$, we can express  the function  $\beta(v(t))$  by $\beta(v(t))=<v,\bar v>=\lambda(t)<\bar v ,\bar v>$. Therefore from  Lemma \ref{mutselecedo.lem-ener1} we deduce that 
\begin{equation}
\lambda'(t)=\frac{1}{K}(\bar \alpha -\alpha(\lambda(t)\bar v +h(t)))\lambda(t). \label{mutselecedo-eq-lambda}
\end{equation}
 Now by using  $\e(h)\to 0$, we deduce that   $h \to 0$ as $t\to \infty$ and   from  \eqref{mutselecedo-eq-lambda} we are reduced to analyse the ODE 
\begin{align*}
\lambda'(t)&=\frac{1}{K}(\bar \alpha -\tilde \alpha(\lambda(t)))\lambda(t)+ \lambda(t)(\tilde \alpha(\lambda(t)\bar v)-  \alpha(\lambda(t)\bar v +h(t)))\\
&=\frac{1}{K}(\bar \alpha -\tilde \alpha(\lambda(t)))\lambda(t)+ \lambda(t) o(1)
\end{align*}
where $|o(1)|=|\tilde \alpha(\lambda(t)\bar v)-  \alpha(\lambda(t)\bar v +h(t))| \le C\sqrt{\e(h)}\to 0$ as $t\to \infty$. As before we can conclude that $\lambda(t)\to 1$ and $v$ converges  to $\bar v$.
 \fdem

\dem{Proof of   Claim \ref{mutselecedo-cla-energy}:} 
Since $\e(h(t))>0$ for all $t$, as in Lemma \ref{mutselecedo.lem-liap} we have
\begin{equation}
\frac{d}{dt}\log\left[\frac{\e(h)}{\left(\beta(v)\right)^2}\right]=  -\frac{1}{\e(h)}\sum_{i,j =1}^{N}\mu_{ij} \bar v_i \bar v_j\left( \frac{h_j}{\bar v_j}-\frac{h_j}{\bar v_j}\right)^2. \label{mutselecedo-cla-dF}
\end{equation} 
Thus the function $\tilde F:= \log\left[\frac{\e(h)}{\left(\beta(v)\right)^2}\right]$ is a decreasing smooth function.

First we observe that the claim is proved if there exists a sequence $(t_n)_{n\in\N}$ converging to infinity so that 
$ \e(h(t_n))\to 0$. Indeed, assume such sequence exists and let  $(s_k)_{k\in \N}$ be a sequence converging to $+\infty$.  Then  there exists $k_0$ and a subsequence $(t_{n_k})_{k\in \N}$ of  $(t_n)_{n\in\N}$ so that  for all $k\ge k_0$,  we have $s_k \ge t_{n_k}$. Therefore from the monotonicity of $\tilde F$ we have for all $k\ge k_0$
$$\log\left[\frac{\e(h(s_k))}{\left(\beta(v(s_k))\right)^2}\right]\le \log\left[\frac{\e(h(t_{n_k}))}{\left(\beta(v(t_{n_k}))\right)^2}\right].$$
By letting $k$ to infinity in the above inequality, we deduce that 
$$\lim_{k\to \infty}\log\left[\frac{\e(h(s_k))}{\left(\beta(v(s_k))\right)^2}\right]=-\infty,$$
which implies that $\e(h(s_k))\to 0, $ since by Lemma \ref{mutselecedo-lem-esti}  $(\beta(v(t_k)))_{k\in \N}$ is uniformly bounded.
The sequence $(s_k)_{k \in \N}$ being chosen arbitrarily this implies that $\e(h(t)) \to 0 $ as $t\to +\infty$.  
\medskip

Let us now prove that such sequence $(t_n)_{n\in \N}$ exists.
We  argue  by contradiction and assume that  $\inf_{t\in \R^+}\e(h(t))>0$. Therefore
from the monotonicity and the smoothness  of $\tilde F$  we deduce that there is  $c_0\in \R$ so that
  $$\tilde F(h(t))\to c_0 \quad \text{ and }\quad \frac{d}{dt}\tilde F(h(t))\to 0 \quad\text{ as }\quad t \to +\infty.$$

Thus by Lemma \ref{mutselecedo-lem-esti} and \eqref{mutselecedo-cla-dF} it follows that 
 
\begin{equation} \lim_{t\to \infty} \sum_{i,j=1}^{N}\mu_{ij}\bar v_j\bar v_i\left(\frac{h_j}{\bar v_j}-\frac{h_i}{\bar v_i}\right)^2=0.\label{mutselecedo-eq-lim}
\end{equation}

From the \textit{ a priori} estimates of Lemma \ref{mutselecedo-lem-esti}, there exists a sequence $t_n\to \infty$ so that for all $i$ $h_i(t_n)\to \tilde h_i$.
Passing to the limit along this sequence in the equation \eqref{mutselecedo-eq-lim} it yields  

$$0=\sum_{i,j=1}^{N}\mu_{ij}\bar v_j\bar v_i\left(\frac{\tilde h_j}{\bar v_j}-\frac{\tilde h_i}{\bar v_i}\right)^2.$$ 

By using the irreducibility assumption on the nonnegative matrix  $\mu_{ij}$  and the positivity of the quantities $\bar v_i$, one can deduce from the above equality that  we must  have for all $i$ and $j$ $$\frac{\tilde h_j}{\bar v_j}=\frac{\tilde h_i}{\bar v_i}.$$ 
Thus if we set $\lambda:=\frac{\tilde h_1}{\bar v_1}$ we have  $ \tilde h=\lambda\bar v$.  So by using that $<h,\bar v>=0$ for all time  it follows that $\lambda=0$. Hence, we get the contradiction $$0<\inf_{t \in \R}\e(h(t))\le \e(h(t_n))\to 0.$$   

\fdem

\section{A case of interest}\label{mutselecedo-section-caseofinterest}

In this section we analyse more precisely the dynamics of the solution $v$ of \eqref{mutselecedo-eq-sysgen} when the  interactions $\Psi_i$ take the form $\Psi_i(v):=\sum_{j=1}^Nr_jv_j$.  We prove the Theorem \ref{mutselecedo-thm2} that we recall below.

\begin{theorem}\label{mutselecedo-thm-coi}
Assume that the interactions $\Psi_i$ take the form $\Psi_i(v)=\sum_{j=1}^Nr_jv_j$ then for any nonnegative  initial datum $v(0)$ not identically  zero, the  solution $v(t)$ of \eqref{mutselecedo-eq-syspart} converges exponentially fast to its unique equilibrium. That is to say there exists two positive constants $C_1$ and $C_2$ so that 
$$\|v-\bar v\|_{\infty}\le C_1e^{-C_2t}.$$ 
\end{theorem}

Before proving  this Theorem we establish two auxiliary results that for convenience  we present in two separate subsections.   
We prove Theorem \ref{mutselecedo-thm2} at the end of this  section.

 \subsection{Study of the evolution of the total population} \label{mutselecedo-sect-totpop}~\\
Let us denote $\n(t)=\sum_{i=1}^{N}v_i$ the total population.
A straightforward computation shows that for the  interactions $\Psi_i$ of the form  $\Psi_i(v):=\sum_{j=1}^Nr_jv_j$ we see that $\n(t)$ satisfies the equation :

 \begin{equation}\label{mutselecedo-eq-totpop}
 \frac{d\n}{dt}= \alpha(t)\left(1-\frac{\n}{K}\right),
 \end{equation}
 
 which written with the new variable  $P(t)=K-\n(t)$ takes the form      
  \begin{equation}\label{mutselecedo-eq-totpop2}
 \frac{dP}{dt}=-\frac{ \alpha(t)}{K}P.
 \end{equation}
 The dynamic of the above equation is strongly related to the behaviour of $\alpha(t)$ and we have 
 \begin{equation}
 |P(t)|=|P_0|e^{-\frac{1}{K}\int_0^t\alpha(s)\,ds}.\label{mutselecedo-eq-totpop3}
 \end{equation}

\begin{lemma}\label{mutselecedo-lem-totpop}
For any nonnegative  initial datum $v(0)$ not identically  zero, $\n(t)$ converges exponentially fast  toward its unique equilibrium $K$. Moreover $\n$ satisfies  identically 
\begin{equation}
\min\{\n_{min}(t),\n_{max}(t)\}\le \n(t)\le \max\{\n_{min}(t),\n_{max}(t)\},\label{mutselecedo-eq-totpop4}
\end{equation} 
where $\n_{min}$ and $\n_{max}$ are   the solutions of the logistic equations: 
 \begin{align*}
&\frac{du}{dt}=\xi^{\pm} u \left(1 -\frac{u}{K}\right)\\ 
&u(0)=\n(0)
\end{align*}
with respectively  $\xi^-=\min\{r_1,\ldots,r_N\}$ and $\xi^+=\max \{r_1,\ldots,r_N\}$ 
\end{lemma}

\dem{Proof:}

Assume for the moment that \eqref{mutselecedo-eq-totpop4} holds then the convergence exponentially fast to $K$ is a straightforward consequence of \eqref{mutselecedo-eq-totpop3}. Indeed  by  \eqref{mutselecedo-eq-totpop4}   we deduce that   $ R_{min}\min( \n_{min},\n_{max}) \le \alpha(s)\le R_{max} \min( \n_{min},\n_{max}) $ where   $R_{\max}:=\max\{r_1,\ldots,r_N\}$ and $R_{\min}:=\min\{r_1,\ldots,r_N\}$. Therefore $\alpha(s)>C_0,$ since $\n_{min}$ and $\n_{max}$ converge to $K$.

To obtain   \eqref{mutselecedo-eq-totpop4} we investigate the following three cases, $\n(0)=K$, $\n(0)<K$ and $\n(0)>K$.

In the first case  $\n(0)=K$, we see that \eqref{mutselecedo-eq-totpop4} holds true trivially, since  $\n_{max}=\n_{\min}=\n(t)\equiv K$ for all $t$. Let us now investigate the two other situations. The argumentation in both situation being similar we expose only the case $\n(0)<K$.

In this situation, $\n_{min}$ and $\n_{max}$ being the solutions of  logistic equations they are  increasing functions. Moreover wa have $\n_{min}(t)\le \n_{max}(t)<K$ for all $t$. 
On another hand, by continuity of $\n(t)$, there exists also $t_1>0$ so that $\n(t)<K$ on $[0,t_1)$. Furthermore  on $(0,t_1)$ we can see that $\n(t)$ satisfies the following differential inequalities
\begin{align}
&\frac{d\n}{dt}\ge R_{min}\n \left(1 -\frac{\n}{K}\right)\\ 
&\frac{d\n}{dt}\le  R_{max}\n \left(1 -\frac{\n}{K}\right).
\end{align}

From the above differentials inequalities, by comparing $\n, \n_{min}$ and $\n_{max}$ via the Cauchy Lipschitz Theorem, we obtain   $\n_{\min}\le \n\le \n_{max}$ for all  $t \in [0,t_1)$. Since $\n_{max}(t_1)<K$, we can bootstrap the above argument and  show that \eqref{mutselecedo-eq-totpop4} holds true for all $t$.

 \fdem

\subsection{A useful functional inequality}~\\
Next we establish a  useful functional inequality satisfied by vectors  $h \in \bar v^{\perp}$ where $\bar v^{\perp}$ denotes the linear subspace of $\R^N$ orthogonal to $\bar v$.

\begin{lemma} \label{mutselecedo-lem-fcineq}
There exists $C_1>0$ so that for all $h\in \bar v^{\perp}$ 
$$C_1\e(h)\le\sum_{i,j =1}^{N}\mu_{ij} \bar v_i \bar v_j\left( \frac{h_j}{\bar v_j}-\frac{h_j}{\bar v_j}\right)^2.$$
Moreover $C_1=\lambda_2$ where $\lambda_2$ is the minimal eigenvalue strictly positive of the linear eigenvalue problem
$$ h_i\sum_{j=1}^N\mu_{ij}\frac{\bar v_j}{\bar v_i}-\sum_{j=1}^N\mu_{ij}h_j=\lambda h_i.$$
\end{lemma}
\dem{Proof :}
Let $\I$ be the following Rayleigh quotient 
$$\I(h):=\frac{1}{\e(h)}\sum_{i,j =1}^{N}\mu_{ij} \bar v_i \bar v_j\left( \frac{h_j}{\bar v_j}-\frac{h_j}{\bar v_j}\right)^2.$$
Observe that the first part of the Lemma is proved if we show that   
\begin{equation}
\inf_{h\in \bar v^{\perp}}\I(h)>0,\label{mutselecedo-eq-energy-expdecay}
\end{equation} 
or equivalently  
\begin{equation}
\inf_{h\in \bar v^{\perp}, \e(h)=1}\I(h)=\inf_{h\in \bar v^{\perp}}\I(h)>0,\label{mutselecedo-eq-energy-expdecay2}
\end{equation} 
since for all real $\mu$, $\I(h)=\I(\mu h)$.

 To obtain \eqref{mutselecedo-eq-energy-expdecay}, we argue by contradiction and assume  that  $\inf_{h\in \bar v^{\perp}}\I(h)=0$.  By \eqref{mutselecedo-eq-energy-expdecay2} we can take $(h_n)_{n\in \N}$  a minimising sequence so  that  $h_n\in\bar v^{\perp},  \e(h_n)=1$ for all $n$. Since $\{x\in \R^{N}\, | x\in \bar v^{\perp}, \e(x)=1\}$ is a closed bounded set, $(h_n)_{n\in \N}$ is uniformly bounded and we can extract a subsequence $(h_{n_k})_{k\in \N}$ which converges to $\bar h \in \bar v^{\perp}$ with $\e(\bar h) =1$. 
 Passing to the limit along this subsequence, we  see that $\I(\bar h)=0$  which combined with $\bar h \in \bar v^{\perp}$ implies that $\bar h =0$. Thus we get the   contradiction  $0=\e(\bar h)=1$.  
Hence \eqref{mutselecedo-eq-energy-expdecay} holds true.

Let $C_1$ be defined by $C_1=\inf_{h\in \bar v^{\perp}, \e(h)=1}\I(h)$, let us try to compute $C_1$. By construction $C_1$  is the result of a constrained minimisation problem. Therefore, by  the standard  optimization Theory \cite{Rockafellar1997},   the minimizers must satisfy the following  Euler-Lagrange equations  
 \begin{align}
 &(D -\tilde M)h=\lambda h,\\
&\e(h)=1,\\
&<h,\bar v>=0,
 \end{align}
    where $\lambda \in \R$ is a Lagrange Multiplier  to be determined and $\tilde M $ and $D$ are the following matrices
    $$ \tilde M:= \begin{pmatrix}  \mu_{11}  & &\mu_{ij} \\
     & \ddots &  \\
   \mu_{ij} &  &  \mu_{NN}\\
\end{pmatrix}, \quad D:= \begin{pmatrix} \frac{1}{\bar v_1}\sum_{j=1}^N\mu_{1j}\bar v_j  & &0 \\
     & \ddots &  \\
   0 &  & \frac{1}{\bar v_N}\sum_{j=1}^N\mu_{Nj}\bar v_j \\
\end{pmatrix}. $$

By taking $k>0$ large enough, the irreducible matrix   $\tilde M -D + k Id$   becomes  non-negative and so $\tilde M -D+kId$ is a symmetric positive definite matrix. As a consequence, the eigenvectors  associated to the eigenvalue $\nu_i$ of $\tilde M -D + k Id$  form an  a orthogonal basis of $\R^N$, i.e.  $\R^N=\bigoplus^{\perp}E_{\nu_i},$ where $E_{\nu_i}$ denotes  the  eigenspace associated to the eigenvalue $\nu_i$.

Since  $\tilde M -D + k Id$ is non negative and irreducible  by the Perron-Frobenius Theorem there exists a unique  eigenvalue, says $\nu_1$ associated with a positive eigenvector $\phi_1$. Furthermore $\nu_1$ is the largest eigenvalue and is algebraically simple.  By a direct computation  one can see that  $(\tilde M -D +kId)\bar v =k\bar v$, thus  $\nu_1=k$ and $\phi_1=\gamma\bar v$. 
From the properties of $\nu_1$ we have $E_{\nu_{1}}=lin(\bar v)$, for all $i\neq 1, \nu_i<\nu_1=k$  and  $\bar v^{\perp}=\bigoplus^{\perp}_{i\neq 1}E_{\nu_i}$.

By construction, the $ \lambda_i:=\nu_i-k \le 0$ are the eigenvalues of $\tilde M -D$ and we can see that for all $h \in \bar v^{\perp}, \e(h)=1$ we have 
$$\I(h)=\transposee{h}(\tilde M -D)h \ge  \min_{i\neq 1}\{-\lambda_i\}.$$

Let $\lambda_2<0$ be second  largest eigenvalue of $\tilde M -D$  and $\phi_2$ an associated eigenvector. 
 By normalising $\phi_2$ properly and since $\phi_2\in \bar v^{\perp}$  a straightforward computation shows that $$\I(\phi_2)=\transposee{\phi_2}(D-\tilde M)\phi_2=-\lambda_2= \min_{i\neq 1}\{-\lambda_i\}.$$
   
Hence, $\min_{h\in \bar v^{\perp}} \I(h)=-\lambda_2.$ 
     
\fdem

\subsection{Asymptotic behaviour of the solution}~\\
We are now in position to obtain the exponential rate of convergence for any solution $v$ of \eqref{mutselecedo-eq-syspart}.

\dem{Proof of Theorem \ref{mutselecedo-thm2}:}

First we claim that 
\begin{equation}\label{mutselecedo-eq-norm-steady}
\sum_{i=1}^N\bar v_i=K.
\end{equation}
 Indeed from Lemma \ref{mutselecedo-lem-steady} we deduce that  $\bar v= \frac{K\lambda_p}{\sum_{i=1}^{N}r_i(v_p)_i}v_p$ where $v_p$ is the positive normalised eigenvector associated to the principal eigenvalue $\lambda_p$ of the matrix $R+M$. 
A straightforward computation shows that 
$$\sum_{i=1}^{N}r_i(v_p)_i=\sum_{i=1}^{N}((R+M)v_p)_i=\lambda_p\sum_{i=1}^{N}(v_p)_i.$$
Since $\lambda_p>0$ it follows  that 
 $$ \frac{\sum_{i=1}^{N}(v_p)_i}{\sum_{i}^{N}r_i(v_p)_i}=\frac{1}{\lambda_p}.$$
Thus $$ \frac{K\lambda_p}{\sum_{i}^{N}r_i(v_p)_i}\sum_{i=1}^{N}(v_p)_i=K.$$ 
 
Now recall that in the proof of Lemma \ref{mutselecedo-lem-asb} from the orthogonal decomposition    $v(t)=\lambda(t)\bar v +h(t)$ we have 
$$\frac{d\e(h)}{dt}=  -\sum_{i,j =1}^{N}\mu_{ij} \bar v_i \bar v_j\left( \frac{h_j}{\bar v_j}-\frac{h_j}{\bar v_j}\right)^2+\frac{2}{K}(\bar \alpha -\alpha(t))\e(h).$$
By  Lemma \ref{mutselecedo-lem-fcineq}   we obtain   
$$\frac{d\e(h)}{dt}\le -C_1\e(h)+\frac{2}{K}(\bar \alpha -\alpha(t))\e(h).$$
Now  arguing as in the proof of Lemma \ref{mutselecedo.lem-liap}, we end up with

  $$\frac{d}{dt}\log\left(\frac{\e(h)}{\beta(v)^2}\right)\le -C_1.$$
 Therefore, thanks to   Lemma \ref{mutselecedo-lem-esti} we deduce that  
 \begin{equation}\label{mutselecedo-eq-asymp1}
 \e(h)\le\left(\frac{\e(h (0))}{\beta(v(0))^2}\right)e^{-C_1t}\beta^{2}(v) \le C_2e^{-C_1t}. 
 \end{equation}

Finally we get the exponential rate of convergence, by observing that $$\n(t)=<v,1>=\lambda(t)<\bar v,1>+<h,1>,$$ 
which thanks to \eqref{mutselecedo-eq-totpop3} implies  that
$$\lambda(t)=\frac{K}{<\bar v,1>} +\frac{<h,1>}{<\bar v,1>}+\frac{|\n(0)-K|e^{-\frac{1}{K}\int_0^t\alpha(s)\,ds}}{<\bar v,1>}.$$ 
By Lemma \ref{mutselecedo-lem-totpop} and the estimates   \eqref{mutselecedo-eq-asymp1}, \eqref{mutselecedo-eq-norm-steady},  ,using  standard norm estimates, we have 
\begin{align*}
 \|v-\bar v\|_{\infty}&\le \left|\frac{<h,1>}{K}+\frac{|\n(0)-K|e^{-\frac{1}{K}\int_0^t\alpha(s)\,ds}}{K}\right| \|\bar v\|_{\infty}+N \e(h)\\
 &\le \frac{|\n(0)-K|\|\bar v\|_{\infty}}{K}e^{-\frac{1}{K}\int_0^t\alpha(s)\,ds} +\left(1+\frac{\|\bar v\|_{\infty}}{K}\right) \tilde C_2e^{-C_1t}.
  \end{align*}

\fdem

\section{The general case:The stationary solution}\label{mutselecedo-section-generalcase}
In this section we investigate the existence of a positive stationary solution of \eqref{mutselecedo-eq-sysgenred} under the additional condition \eqref{mutselecedo-hyp3} on the matrices $R$ and $M$ that we recall below,
 $$ \forall\; i\quad \sum_{j=1}^{N}\mu_{ij}\le \frac{r_i}{2}.$$

This assumption has for consequence that the matrix $R+M$ is positive definite. Indeed, we have
\begin{align*}
 \transposee{h}(R+M)h=\sum_{i=1}^{N}r_i h_i^2 - \frac{1}{2}\sum_{i,j=1}^{N}\mu_{ij}(h_i-h_j)^2 &\ge 2\sum_{i,j=1}^{N}\mu_{ij}h_i^2  - \frac{1}{2}\sum_{i,j=1}^{N}\mu_{ij}(h_i-h_j)^2\\
&\ge  \sum_{i,j=1}^N(\mu_{ij}\left(h_i^2+h_j^2- \frac{1}{2}(h_i-h_j)^2\right)\\
&\ge\frac{1}{2}\left( \sum_{i,j=1}^{N}\mu_{ij}(h_i+h_j)^2\right).
 \end{align*}     
Thus  $kern(R+M)=\{0\}$  and the matrix  $R+M$ is invertible. Moreover from the last inequality we see that there exists
positive constants $c_0$ and $C_0$ so that
\begin{equation}
 c_0<u,u>_{R+M}\le \e(u) \le C_0<u,u>_{R+M}, \label{mutselecedo-eq-equiv-scalarp}
 \end{equation}
 where $<u,u>_{R+M}:=\transposee{u}(R+M)u$.
 
Let  $\Xi(v)$  be the diagonal matrix defined by 
$$(\Xi(v))_{ij}=\delta_{ij}\Psi_i(v).$$
With this notation, a positive stationary solution of \eqref{mutselecedo-eq-sysgenred} is then  a non negative solution of the following problem:

\begin{equation}\label{mutselecedo.eq-gen-steady}
(R+M)v=\Xi(v)v
\end{equation}

Note that when $\Xi(v)$ can be written  as $\Xi(v)=\alpha(v)Id$, the construction of a positive solution has already been made in Section \ref{mutselecedo-section-specialcase}. So in the later, we will assume that  $\Xi(v)$ cannot be written as $\Xi(v)=\alpha(v)Id$.
It is worth mentioning that in this situation the method used in Section \ref{mutselecedo-section-specialcase} does not  work and we have to use another strategy.

Let $T$ be the following map

$$\begin{array}{rcl}
T: \R^{N}&\to&\R^{N}\\
\\v&\mapsto&Tv:=(R+M)^{-1}[\Xi(v)v]\end{array}.$$
Since $R+M$ is invertible, $T$ is well defined and  one can easily check that a positive solution of \eqref{mutselecedo.eq-gen-steady} is a positive fixed point of the map $T$.
To check that $T$ has a positive fixed point we use a degree argument.

 Let $\Psi_i(v)^s$ and $\Xi(v)^s$ defined by  
$$
\Psi_{i}^s:=s\Psi_i+(1-s)\Psi_1,\qquad (\Xi^s(v))_{ij}:=\delta_{ij}\Psi^s_i(v),
$$
we consider the homotopy  $H \in C([0,1]\times \R^{N}, \R^{N})$ defined by
$$\begin{array}{ccl}
H: [0,1]\times\R^{N}&\to&\R^{N}\\
(s,v)&\mapsto&H(s,v):=(R+M)^{-1}[\Psi^s(v)v].\end{array}$$
One can see that $H(1,.)=T$ and $H(0,.)=T_0$ where $T_0$ corresponds to the map 
$$\begin{array}{rcl}
T_0: \R^{N}&\to&\R^{N}\\
\\v&\mapsto&T_0v:=\Psi_1(v)(R+M)^{-1}v.\end{array}$$

Note that there exists a unique positive fixed point to $T_0$ which can be constructed by arguing as in Section \ref{mutselecedo-section-specialcase}.

The next step in this degree argument is to obtain  for all $s$, a good \text{ a priori } estimates on the fixed point of the map $H(s,v)$, i.e.  a good estimate on the  positive solutions of the following problem:
\begin{equation}
(R+M)V=\Psi^s(V)V. \label{mutselecedo-eq-homotop-fp}
\end{equation}

 In this direction we show the following:

\begin{lemma}\label{mutselecedo-lem-esti1}
Let $V$ be a non negative solution of  \eqref{mutselecedo-eq-homotop-fp}. Then either $V\equiv 0$ or $V>0$ and there exists $\bar c_1$ and $\bar C_1$ independent of $s$  so that $$\bar c_1\le \sum_{i}V_i\le \bar C_1. $$
\end{lemma}

\dem{Proof:}
To obtain that $V$ is either positive or $V=0$   we can argue as in  the proof of Lemma \ref{mutselecedo-lem1}.  So  assume  that there exists a $i_0$ so that $V_{i_0}=0$.  By construction  $V_{i_0}$ is a minimum of the $V_i$ and from the equation satisfied by $V_{i_0}$ we get 
$$ 0\le \sum_{j=1}^N\mu_{i_0j}(V_j-V_{i_0})=0.$$
Therefore $V_j=V_{i_0}$ for all $j$ where $\mu_{i_{0}j}\neq 0$. Since $M$ is irreducible there exists $j\neq i_0$ so that $\mu_{i_0 j}\neq 0$. Let $\rho:=\{k\,|\, V_k=0\}$ then $i_0$ and all $j$ so that $\mu_{i_0j}\neq 0$ belongs to $\rho$. In the previous argument, by replacing $i_0$ by any $k\in \rho$, we see that all $j$ so that $\mu_{kj}\neq 0$ belongs to the set $\rho$.  By iterating enough times the above argument and using the irreducibility of the  matrix $M$ we can see that $\rho=\{1,\ldots, N\}$ so $V_i=V_{i_0}=0$ for all $i$.  
Therefore a nonnegative solution of  \eqref{mutselecedo-eq-homotop-fp} is either $V\equiv 0$ or $V>0$. 

Now let us assume that $V>0$.
Recall that by \eqref{mutselecedo-eq-equiv-scalarp} there exists positive constants $c_0$ and $C_0$ so that for all $u\in \R^N$
$$ c_0<u,u>_{R+M}\le \e(u) \le C_0<u,u>_{R+M}.$$
 So for a solution $V$ of \eqref{mutselecedo-eq-homotop-fp} one has
\begin{align*}
<V,V>_{R+M}=\sum_{i}\Psi_i^s(V)V_i^2\ge &\left(s\min_{i\in \{2,\ldots , N\}}\Psi_i(V)+\Psi_1(V)\right)<V,V>\\
\ge &c_0 \Psi_1(V) <V,V>_{R+M},
\end{align*}
and we also  get 
$$
<V,V>_{R+M}\le C_0\left(s\max_{i\in \{2,\ldots , N\}}\Psi_i(V)+ \Psi_1(V)\right) <V,V>_{R+M}.
$$

Therefore we have 
\begin{align}
&\Psi_1(V)\le \frac{1}{c_0}, \label{mutselecedo-eq-esti1Psi1}\\
&\frac{1}{C_0}\le s\max_{i\in \{2,\ldots , N\}}\Psi_i(V)+ \Psi_1(V)\label{mutselecedo-eq-estilip}.
\end{align}

Now thanks to  the assumptions \eqref{mutselecedo-hyp1}--\eqref{mutselecedo-hyp2} made on the functions $\Psi_i$,  there exists $R_1, c_1, k_1$ and $N$ positive constants $\kappa_i$ so that :

\begin{align}
&\text{ For all }\;  x\in \R^{N,+}\setminus Q_{R_1}(0),\quad c_{1}\left(\sum_{i=1}^Nx_i\right)^{k_{1}}  \le \Psi_1(x), \label{mutselecedo-eq-esti2Psi1}\\ 
&\text{ For all } i,\; \text{and for all} \;  x\in Q_{R_1}(0),\quad \Psi_i(x)\le \kappa_i  \sum_{j=1}^N|x_j| . \label{mutselecedo-eq-esti3Psi1} 
\end{align}

By combining \eqref{mutselecedo-eq-esti1Psi1},\eqref{mutselecedo-eq-esti2Psi1},  \eqref{mutselecedo-eq-estilip} and \eqref{mutselecedo-eq-esti3Psi1} we deduce that  
\begin{align*}
& \n=\sum_{j=1}^{N}V_j\le \sup\left\{\left(\frac{1}{c_0c_1}\right)^{\frac{1}{k_1}},R_1\right\},\\
&\n=\sum_{j=1}^{N}V_j\ge  \min \left\{R_1, \frac{1}{C_0(\kappa_1+\sup_{i}\kappa_i)} \right\}.
\end{align*}

\fdem
\subsection{Computation of the degree}~\\

We are now in position to prove the existence of a positive solution to the equation \eqref{mutselecedo.eq-gen-steady} by means of the computation of the topological degree of $T-id$ on a well chosen set $\O\subset\R^{N,+}$.
Now  we take  two positive constants $c_2$ and $C_2$ so that $c_2<\bar c_1$ and $C_2>\bar C_1$ where $\bar c_1$ and $\bar C_1$ are the constants obtained in  Lemma \ref{mutselecedo-lem-esti1}.  Let $\O$ be the following open set
$$\O:=\left\{v\in\R^{N,+}\,|\,   c_2\le \sum_{i=1}^{N}v_i \le C_2 \right\}$$  
and let us compute $deg(T-Id, \O,0)$.

By Lemma \ref{mutselecedo-lem-esti1}  for all $s \in [0,1]$  in $ H(s,v)-v\neq 0$ on $\partial \O$. Therefore using that  $H(.,.)$ is an homotopy,    we conclude that  $deg(T-Id, \O,0)=deg(H(1,.)-Id,\O,0)=deg(H(0,.)-Id,\O,0).$
By construction,  from Section \ref{mutselecedo-ss-equilibria}, one can check that $deg(H(0,.)-Id,\O,0)\neq 0 $ since  the map $T_0$ has an unique stable and positive fixed point.
Thus $deg(T-Id, \O,0)\neq 0$ which shows that $T$ has a fixed point in $\O$.
\fdem

\section{The general case: Asymptotic Behaviour}\label{mutselecedo-section-genasb}
In this section we prove Theorem \ref{mutselecedo-thm4}.  That is to say, under the extra assumption \eqref{mutselecedo-hyp3} we analyse the asymptotic behaviour of the solution $v(t)$ when for  all $i$ the interaction $\Psi_i$ can be expressed like: $\Psi_i(v)=\alpha(v)+\eps \psi_i(v)$ with  $\psi_i$  uniformly bounded.
To obtain the asymptotic behaviour in this case,  we follow the strategy  developed  in Subsection \ref{mutselecedo-ss-asb}. Namely, we start by showing some \textit{a priori} estimates on the solution,  then we analyse the convergence by means of a Lyapunov functional. For convenience we dedicate a subsection to each essential part of the proof.  
 \smallskip
 
\subsection{A priori estimate }~\\ 
We start by establishing some useful differential inequalities. Namely we show that 
\begin{lemma}\label{mutselecedo-lem-diffineq}
Assume that $r_i,(\mu_{ij}), \Psi_i$ satisfies \eqref{mutselecedo-hyp1}, \eqref{mutselecedo-hyp2} and  \eqref{mutselecedo-hyp3}. Assume further that $\Psi_i(v)=\alpha(v)+\eps \psi_i(v)$ with $\psi_i$  uniformly bounded.  Then there exists $\eps_0>0$ so that for all $0\le \eps\le\eps_0$ there exists $\o^+\in \R^N,$ $\o^+$ positive and  a positive real $\gamma$ so that 
\begin{itemize}
\item[(i)] 
\begin{align*}
&\frac{d\sum_{i=1}^N \bar\o^+_iv_i}{dt}\le    \frac{1}{K} (\alpha( \bar\o^+)-\alpha(v)) \sum_{j=1}^N \bar\o^+_iv_i\\
&\frac{d\sum_{i=1}^N \bar \o^+_iv_i}{dt}\ge    \frac{1}{K} (\alpha(\gamma \bar\o^+)-\alpha(v)) \sum_{j=1}^N \bar\o^+_iv_i
\end{align*}
\item[(ii)]
\begin{align*}
&\frac{d\e(v)}{dt}\le -\sum_{i,j =1}^{N}\mu_{ij}  \bar\o^+_i  \bar\o^+_j\left( \frac{v_i}{ \bar\o^+_i}-\frac{v_j}{ \bar\o^+_j}\right)^2 +\frac{2}{K}(\alpha(\o^+)-\alpha(v)) \e(v)\\
 &\frac{d\e(v)}{dt}\ge -\sum_{i,j =1}^{N}\mu_{ij}  \bar\o^+_i  \bar \o^+_j\left( \frac{v_i}{ \bar\o^+_i}-\frac{v_j}{ \bar\o^+_j}\right)^2 +\frac{2}{K}(\alpha(\gamma \bar\o^+)-\alpha(v)) \e(v)
 \end{align*}
 
 \item[(iii)]
 \begin{align*}
 &\frac{d}{dt}\log\left[\frac{\e(v)}{(\sum_{i=1}^N  \bar\o^+_iv_i)^2}\right]\le -\frac{1}{\e(v)}\sum_{i,j =1}^{N}\mu_{ij}  \bar\o^+_i  \bar\o^+_j\left( \frac{v_i}{ \bar\o^+_i}-\frac{v_j}{ \bar\o^+_j}\right)^2 +\frac{2}{K}(\alpha( \bar\o^+)-\alpha( \gamma \bar\o^+))\\
  &\frac{d}{dt}\log\left[\frac{\e(v)}{(\sum_{i=1}^N  \bar\o^+_iv_i)^2}\right]\ge -\frac{1}{\e(v)}\sum_{i,j =1}^{N}\mu_{ij}  \bar\o^+_i \bar \o^+_j\left( \frac{v_i}{ \bar\o^+_i}-\frac{v_j}{ \bar\o^+_j}\right)^2 +\frac{2}{K}(\alpha(\gamma \bar\o^+)-\alpha(  \bar\o^+))
\end{align*}
\end{itemize}
\end{lemma}

\dem{Proof:}
First, we observe that $(iii)$ can be straightforwardly obtained by combining $(i)$ and $(ii)$. So we deal only with $(i)$ and $(ii)$.

Let us denote $\sigma:=\eps\|\psi\|_{\infty}$. Then since $v$ is positive, from \eqref{mutselecedo-eq-sysgen} it follows that 
\begin{align*}
\frac{dv_i}{dt}&\le (r_i+\sigma -\alpha(v))v_i + \sum_{j=1}^N\mu_{ij}(v_j -v_i),\\
\frac{dv_i}{dt}&\ge(r_i-\sigma -\alpha(v))v_i + \sum_{j=1}^N\mu_{ij}(v_j -v_i).
\end{align*}

Let $\bar \o^+$ and $\bar \o^-$ be the stationary solutions of the corresponding equations   
\begin{align*}
&\frac{d\o_i^+}{dt}=(r_i+\sigma -\alpha(\o^+))\o^+_i + \sum_{j=1}^N\mu_{ij}(\o^+_j -\o^+_i),\\
&\frac{d\o_i^-}{dt}=(r_i-\sigma -\alpha(\o^-))\o^-_i + \sum_{j=1}^N\mu_{ij}(\o^-_j -\o^-_i).
\end{align*}
Now, let us fix  $\eps$ small enough so that $\bar \o^{\pm}$ exists.
Then by arguing as in the proof of Lemma \ref{mutselecedo-lem-general-id}, we obtain 
\begin{align*}
&\frac{d}{dt}\sum_{i=1}^N  \bar\o^+_iv_i\le    \frac{1}{K} (\alpha( \bar\o^+)-\alpha(v)) \sum_{j=1}^N \bar\o^+_iv_i,\\
&\frac{d}{dt}\sum_{i=1}^N  \bar\o^-_iv_i\ge    \frac{1}{K} (\alpha( \bar\o^-)-\alpha(v)) \sum_{j=1}^N \bar\o^-_iv_i,\\
&\frac{d\e(v)}{dt}\le -\sum_{i,j =1}^{N}\mu_{ij} \bar \o^+_i  \bar\o^+_j\left( \frac{v_i}{ \bar\o^+_i}-\frac{v_j}{ \bar\o^+_j}\right)^2 +\frac{2}{K}(\alpha( \bar\o^+)-\alpha(v)) \e(v),\\
 &\frac{d\e(v)}{dt}\ge -\sum_{i,j =1}^{N}\mu_{ij} \bar \o^-_i  \bar\o^-_j\left( \frac{v_i}{ \bar\o^-_i}-\frac{v_j}{ \bar\o^-_j}\right)^2 +\frac{2}{K}(\alpha( \bar\o^-)-\alpha(v)) \e(v).
\end{align*}

Note that by Lemma \ref{mutselecedo-lem-steady}, $\bar \o^{\pm}$ are   positive eigenvectors of the matrices $R+M\pm \sigma Id$. Thus $\bar \o^{\rm}$ are  eigenvectors of the matrix $R+M$ associated to the principal eigenvalue of $R+M$. Since $R+M$ is irreducible,  the eigenspace  associated to the principal eigenvalue   is unidimensional. So, we have   $\o^+=\gamma\o^-$ for some positive $\gamma$.
 Hence $(i)$ and $(ii)$ hold true since 
 $$ \sum_{i,j =1}^{N}\mu_{ij} \bar \o^-_i  \bar\o^-_j\left( \frac{v_i}{ \bar\o^-_i}-\frac{v_j}{ \bar\o^-_j}\right)^2=\sum_{i,j =1}^{N}\mu_{ij} \bar \o^+_i  \bar\o^+_j\left( \frac{v_i}{ \bar\o^+_i}-\frac{v_j}{ \bar\o^+_j}\right)^2.$$

\fdem
\medskip

Next, we derive some \textit{ a priori } estimates for the solution $v$ for an interaction $\Psi$ as in Theorem \ref{mutselecedo-thm4}.
We first prove some sharp  \textit{ a priori } estimates for stationary solution $\bar v_\eps$ of \eqref{mutselecedo-eq-sysgen}.

\begin{lemma}\label{mutselecedo-lem-estigen1}
Assume that $r_i,(\mu_{ij}), \Psi_i$ satisfies \eqref{mutselecedo-hyp1}, \eqref{mutselecedo-hyp2} and  \eqref{mutselecedo-hyp3}. Assume further that $\Psi_i(v)=\alpha(v)+\eps \psi_i(v)$ with $\psi_i$  uniformly bounded.  Then  there exists $\bar c_1<\bar C_1$ and $\eps_1$  so that  for all $0\le \eps\le \eps_1$ and for any positive stationary solution $\bar v_\eps$  of \eqref{mutselecedo-eq-sysgen}, we have 
$$\bar c_1\le \sum_{i=1}^{N}\bar v_{\eps,i}<\bar C_1.$$  
Moreover, there exists $r_0$ so that for all $\eps \le \eps_1, \bar v_\eps \in Q_{\sqrt{\eps} r_0}(\bar \o^+)$. 
Furthermore, for any nonnegative initial datum $v_i(0)$  not identically  zero there exists two constants $\bar c_2(v(0)), \bar C_2(v(0))$ so that for all $0\le\eps\le \eps_1$ 
the solution $v_\eps$ satisfies for all $t,$
  $$ \bar c_2\le \beta(v_\eps(t)):=\sum_{i=1}^N\bar v_{\eps,i} v_{\eps,i}(t) \le \bar C_2. $$
\end{lemma}

\dem{Proof:}

Let us first observe that for $\eps\le \eps_0$ by replacing $v_\eps$ by $\bar v_\eps$ in $(i)$ of Lemma \ref{mutselecedo-lem-diffineq}, we get 
\begin{align}
&0\le \frac{1}{K} (\alpha( \bar\o^+)-\alpha(\bar v_\eps)) \sum_{j=1}^N \bar\o^+_i \bar v_{\eps,i},\label{mutselecedo-eq-esti0-eps}\\
&0\ge   \frac{1}{K} (\alpha( \gamma\bar\o^+)-\alpha(\bar v_\eps)) \sum_{j=1}^N \bar\o^+_i\bar v_{\eps,i}.\label{mutselecedo-eq-esti00-eps}
\end{align}

From the proof of Lemma \ref{mutselecedo-lem-steady}, we also see that  
\begin{equation}
 \alpha(\bar \o^+ ) = K(\lambda_p+\sigma),\qquad \alpha(\bar \o^-)=K(\lambda_p -\sigma),  \label{mutselecedo-eq-esti1-eps}
\end{equation}
where $\lambda_p$ is the positive principal eigenvalue of the matrix $R+M$.
By using the monotonicity of $\alpha$,  we deduce that the maps $\sigma\mapsto \o^{\pm}$ are monotone. 
Moreover, we have $0\le \alpha(\bar \o^+) -\alpha(\bar \o^-)\le 2\sigma K$.

To obtain $\bar c_1$ and $\bar C_1$ we argue as follows.
Let us fix $ \eps_1:= \min \{\eps_0, \frac{\lambda_p}{4\|\psi\|_{\infty}}\}$. Then by \eqref{mutselecedo-eq-esti1-eps}, for all $\eps\in [0,\eps_1]$  we have 

\begin{equation}
 \alpha(\bar \o^+ ) \le \alpha(\o^+_{\eps_1})=\frac{5K\lambda_p}{4},\qquad \alpha(\bar \o^-)\ge \alpha(\o^{-}_{\eps_1})=\frac{3K\lambda_p}{4}.  \label{mutselecedo-eq-esti2-eps}
\end{equation}

Now thanks to  the assumptions \eqref{mutselecedo-hyp1}--\eqref{mutselecedo-hyp2} satisfied by $\alpha$,  there exists $R_\alpha, c_\alpha, k_\alpha$ and $\kappa_\alpha$ so that :

\begin{align}
&\text{ For all }\;  x\in \R^{N,+}\setminus Q_{R_\alpha}(0),\quad c_{\alpha}\left(\sum_{i=1}^Nx_i\right)^{k_{\alpha}}  \le \alpha(x), \label{mutselecedo-eq-esti3-eps}\\ 
&\text{ For all } i,\; \text{and for all} \;  x\in Q_{R_\alpha}(0),\quad \alpha(x)\le \kappa_{\alpha}  \sum_{j=1}^N|x_j| . \label{mutselecedo-eq-esti4-eps} 
\end{align} 
 
So  for $\eps \in [0,\eps_1]$,  by combining \eqref{mutselecedo-eq-esti0-eps}, \eqref{mutselecedo-eq-esti00-eps},\eqref{mutselecedo-eq-esti1-eps},\eqref{mutselecedo-eq-esti2-eps},  \eqref{mutselecedo-eq-esti3-eps} and \eqref{mutselecedo-eq-esti4-eps} we achieve  
\begin{align}
& \sum_{j=1}^{N}\bar v_{\eps,j}\le \sup\left\{\left(\frac{5K \lambda_p}{4c_\alpha}\right)^{\frac{1}{k_\alpha}},R_\alpha\right\}=:\bar C_1,\label{mutselecedo-eq-esti5-eps}\\
&\sum_{j=1}^{N}\bar v_{\eps,j}\ge  \min \left\{R_\alpha, \frac{3K\lambda_p}{4\kappa_\alpha} \right\}=:\bar c_1.\label{mutselecedo-eq-esti6-eps}
\end{align}

To obtain a more precise estimate on $\bar v_\eps$, we argue as follows. Let us decompose $\bar v_\eps:=\lambda_\eps \bar \o^++h_\eps$ where $h_\eps$ is orthogonal to $\bar \o^+$. Then by replacing $v_\eps$ by $\bar v_\eps$ in $(iii)$ of Lemma \ref{mutselecedo-lem-diffineq}, and using that $|\alpha(\bar \o^+) -\alpha(\bar \o^-)|\le 2\sigma K $ we get 
\begin{align}
 0&\le -\frac{1}{\e(\bar v_\eps)}\sum_{i,j =1}^{N}\mu_{ij}  \bar\o^+_i  \bar\o^+_j\left( \frac{h_{\eps,i}}{ \bar\o^+_i}-\frac{h_{\eps,j}}{ \bar\o^+_j}\right)^2 +4\sigma,\label{mutselecedo-eq-esti0-mp}\\
0  &\ge -\frac{1}{\e(\bar v_\eps)}\sum_{i,j =1}^{N}\mu_{ij}  \bar\o^+_i \bar \o^+_j\left( \frac{h_{\eps,i}}{ \bar\o^+_i}-\frac{h_{\eps,j}}{ \bar\o^+_j}\right)^2 -4\sigma.\label{mutselecedo-eq-esti1-mp}
\end{align}
From \eqref{mutselecedo-eq-esti5-eps} and using the functional inequality, Lemma \ref{mutselecedo-lem-fcineq}, we deduce that 
\begin{equation}
\e(h_\eps)\le \frac{4\sigma \bar C_1^2}{C(\bar \o^+)}\label{mutselecedo-eq-esti-2-mp}.
\end{equation}
Combining the latter estimate with \eqref{mutselecedo-eq-esti5-eps} and the positivity of $\bar\o^+$ and $\bar v_\eps$, we have the estimate 
$$\lambda_\eps\le \frac{\bar C_1}{\sum_{i=1}^N\bar \o_i}\left(1+2\sqrt{\frac{N\sigma}{C(\bar \o)}}\right), $$  
where $\bar \o$ is the stationary solution with $\eps=0$. Let $R_0:=\bar C_1\left(1+2\sqrt{\frac{N\sigma}{C(\bar \o)}}\right),$ and choose $\eps_1$ smaller if necessary to have $R_0\le 2\bar C_1$. Next consider the set $\bar Q_{R_0}(0)\subset Q_{2\bar C_1}(0)$. From the above estimates on $\lambda_{\eps}$ and $h_\eps$, since $\alpha$ is Lipschitz continuous in $Q_{2\bar C_1}(0)$ there exists $\kappa_0$ so that for all $\eps$, 
$$
|\alpha(\lambda_\eps \bar \o^+) -  \alpha(\lambda_\eps \bar \o^++h_\eps)|\le \kappa_0 \sqrt{N \e(h_\eps)},
$$
which combine with \eqref{mutselecedo-eq-esti-2-mp} enforces
\begin{equation}\label{mutselecedo-eq-esti-3-mp}
|\alpha(\lambda_\eps \bar \o^+) -  \alpha(\lambda_\eps \bar \o^++h_\eps)|\le 2\bar C_1\kappa_0 \sqrt{\frac{N \sigma}{C(\bar \o)}}.
\end{equation}
Thus from  \eqref{mutselecedo-eq-esti0-eps}, \eqref{mutselecedo-eq-esti00-eps} and \eqref{mutselecedo-eq-esti-3-mp} 
we deduce that
\begin{equation}
|\alpha(\bar\o^+)-\alpha(\lambda_\eps\bar \o^+)|\le 2K\sigma + 2\bar C_1\kappa_0 \sqrt{\frac{N \sigma}{C(\bar \o)}}\le C\sqrt{\sigma},\label{mutselecedo-eq-esti-4-mp}
\end{equation}
with $C$ independent of $\eps$. %
 %
 Observe that by construction there exists two positives constants $0<\iota_0<1<\iota_1$ so that for all $\eps \le \eps_1$, there exists $\iota_\eps\in (\iota_0,\iota_1)$ so that  $\bar\o^+= \iota_\eps \bar \o$. 
 Recall now that by assumption $\nabla\alpha>0$, then the real map  $s \mapsto \alpha(s\bar \o)$ is  smooth ($C^1(\R)$) and increasing.  It is therefore an homeomorphism in $\R^+$ and a local diffeomorphism in $\R^+$.  
 So by the Inverse Function  Theorem, we deduce that for all $s,t \in (0,C)$,
  $$|s-t|\le \frac{k}{\iota_0}|\alpha(s\bar \o^+)-\alpha(t\bar \o^+)|, $$   
where \begin{align*}
&C:= sup\left\{\frac{\iota_1}{\iota_0},\frac{\bar C_1}{\sum_{i=1}^N\bar \o_i}\left(1+2\sqrt{\frac{N\sigma}{C(\bar \o)}}\right)\right\}\\
&k:=\frac{1}{\min_{s\in(0,C)}(<\nabla(\alpha(s\bar\o)),\bar \o>)}.   
\end{align*}
In particular, we have 
$$|s-1|\le \frac{k}{\iota_0}|\alpha(s\bar \o^+)-\alpha(\bar \o^+)|, $$
which combined with \eqref{mutselecedo-eq-esti-4-mp} enforces
$$|\lambda_\eps -1|\le  \frac{k}{\iota_0}C\sqrt{\sigma}.$$
Hence  from the latter estimate and \eqref{mutselecedo-eq-esti-2-mp} we have for all $\eps \le \eps_1$, $\bar v_\eps \in Q_{C\sqrt{\eps}}(\bar \o^+)$.


Next, we derive an uniform upper bound for $\beta(v_\eps)$ when $\eps\in [0,\eps_1]$. In the sequel of this proof, for convenience we drop the subscript $\eps$ on $v$.   

First, we observe that by Lemma \ref{mutselecedo-lem-general-id}  we have 
\begin{align*}
 \frac{d}{dt}\left(\sum_{i=1}^{N}v_i\bar v_i \right)&=\frac{1}{K}\sum_{i=1}^N(\Psi_i(\bar v)-\Psi_i(v))v_i\bar v_i,\\
&=\left(\frac{1}{K}(\alpha(\bar v)-\alpha(t))\right) \sum_{i=1}^Nv_i\bar v_i + \frac{\eps}{K}\sum_{i=1}^N(\psi_i(\bar v)-\psi_i(v))v_i\bar v_i.
\end{align*}
Since the functions $\psi_i$ are uniformly bounded, we achieve
$$ \frac{d}{dt}\left(\sum_{i=1}^{N}v_i\bar v_i \right)\le \left(\frac{1}{K}(\alpha(\bar v)+2\eps\|\psi\|_{\infty}-\alpha(t))\right) \sum_{i=1}^Nv_i\bar v_i .$$
By using \eqref{mutselecedo-eq-esti0-eps} and \eqref{mutselecedo-eq-esti2-eps}, it follows
   
$$ \frac{d}{dt}\left(\sum_{i=1}^{N}v_i\bar v_i \right)\le \frac{1}{K}\left[\frac{7K\lambda_p}{4}-\alpha(t)\right] \sum_{i=1}^Nv_i\bar v_i .$$
Again using that  $\alpha$ satisfies the assumptions \eqref{mutselecedo-hyp1}--\eqref{mutselecedo-hyp2},  there exists $R_\alpha, c_\alpha, k_\alpha$ so that for all $x\in \R^{N,+}\setminus Q_{R_\alpha}(0)$,
\begin{equation}
c_{\alpha}\left(\sum_{i=1}^Nx_i\right)^{k_{\alpha}}  \le \alpha(x). \label{mutselecedo-eq-estieps-alpha} 
\end{equation}

Now, let us assume that for some $t>0, v(t) \in  \R^{N,+}\setminus \bar Q_{R_\alpha}(0)$ otherwise the proof is done since we have
$$\beta(v)\le \left(\max_{i\in \{1,\ldots, N\}}\bar v_i\right) \sum_{j=1}^N|v_j| \le \bar C_1R_{\alpha},$$
where $\bar C_1$ is the bound obtained above.
Let   $\Sigma$ be the following set
$$
\Sigma:=\{t\in \R^+\, |\, \n(t)>R_{\alpha}\}.
$$
 From \eqref{mutselecedo-eq-estieps-alpha}, \eqref{mutselecedo-eq-esti5-eps} and using that the  $v_i$ are non negative,   we see that  for all  $ t  \in\Sigma$
 $$ \alpha(t)\ge c_\alpha\n(t)^{k_\alpha} \ge c_\alpha\left(\frac{1}{\bar C_1}\right)^{k_\alpha}\beta(t)^{k_\alpha}. $$
Therefore, with $\tilde c_0:=c_\alpha\left(\frac{1}{\bar C_1}\right)^{k_\alpha}$  we have for  $t\in\Sigma$
 $$\frac{d \beta(t)}{dt}\le \frac{1}{K}\left(\frac{7K\lambda_p}{4} -\tilde c_0\beta^{k_{\alpha}}(t)\right)\beta(t).$$
 Using the logistic character of the above equation, we can check that 
 $$\beta(v(t))\le \sup\left\{\beta(v(0)),\sup_{x\in Q_{R_{\alpha}}(0)}\beta(x),\left(\frac{7K\lambda_p}{4\tilde c_0}\right)^{\frac{1}{k_{\alpha}}}\right\}.$$
 
 Thus, by using \eqref{mutselecedo-eq-esti5-eps},  we achieve for all $\eps \in [0,\eps_1]$
 $$\beta(v_\eps(t))\le \bar C_2:=\sup\left\{\beta(v(0)),\bar C_1R_\alpha,\bar C_1\left(\frac{7K\lambda_p}{4c_\alpha}\right)^{\frac{1}{k_{\alpha}}}\right\}.$$

To obtain the lower bound for $\beta(v_\eps)$ we argue as follows. First let us observe that by taking $\eps_1$ smaller if necessary, since $\bar \o^+ \to \bar \o$ as $\eps \to 0$ and $\bar v_\eps \in Q_{C\sqrt{\eps}}(\bar \o^+)$, there exists a positive constant $c_0$ independent of $\eps \in [0,\eps_1]$ so that for all stationary solution $\bar v_\eps$ we have 
$$ \min_{i\in\{1,\ldots,N\}}\bar v_{\eps,i}\ge c_0.$$ 
Now by \eqref{mutselecedo-eq-esti6-eps}, for all $\eps \in [0,\eps_1]$ we have 
\begin{equation}
\beta(v_\eps)\ge \min_{i\in\{1,\ldots,N\}}\bar v_{\eps,i}\sum_{i=1}^Nv_{\eps,i}\ge c_0\sum_{i=1}^Nv_{\eps,i}. \label{mutselecedo-eq-esti-betaeps}\end{equation}

 Therefore to obtain an uniform lower bound for $\beta(v_\eps)$ it is enough to obtain an uniform lower bound for $\n_\eps:=\sum_{i=1}^N v_{\eps,i}$.
 From \eqref{mutselecedo-eq-sysgen} by summing over all $i$ and by using the definition of $\Psi$ and the boundedness of the $\psi_i$ we deduce that 
$\n_\eps$ satisfies the following inequality
 $$ \frac{d \n_{\eps}}{dt}\ge\frac{1}{K}\left(\alpha(\bar v_\eps)-2\eps\|\psi\|_{\infty}-\alpha(t)\right) \n_{\eps} .$$
  
Thanks to \eqref{mutselecedo-eq-esti0-eps} and \eqref{mutselecedo-eq-esti2-eps}, we have 
$$ \frac{d \n_{\eps}}{dt}\ge\frac{1}{K}\left(\frac{K\lambda_p}{4}-\alpha(t)\right) \n_{\eps} .$$

 By reproducing the argumentation  of the proof of Lemma \ref{mutselecedo-lem1} and by using Remark \ref{mutselecedo-rem-lem1} we can  check that 
 \begin{equation}\label{mutselecedo-eq-esti-Neps}
 \n_{\eps}(t)\ge \min \left\{1,\frac{K\lambda_p}{4\kappa_{1}}, \frac{\n(0)}{2}\right\},
 \end{equation}
where $\kappa_1$ denotes the Lipschitz constant of the function $\alpha$ in the unit cube. 
Hence, by collecting \eqref{mutselecedo-eq-esti-betaeps} -\eqref{mutselecedo-eq-esti-Neps}  we achieve for all $\eps \le \eps_1$ and all $t>0$,
$$ \beta(v_\eps(t))\ge \frac{\bar c_1}{N}\min \left\{1,\frac{K\lambda_p}{4\kappa_{1}}, \frac{\sum_{j=1}^Nv_i(0)}{2}\right\}=:\bar c_2. $$

\fdem
\medskip

\begin{remark} Note that from the above argumentation, using the Logistic character of the equations,  we can  get that  for all $\eps \le \eps_1$ and all  initial data $v(0)\ge 0$, there exists $t_0$ so that for all $t \ge t_0$ we have 
$$\frac{1}{2}\frac{\bar c_1}{N}\min \left\{1,\frac{K\lambda_p}{4\kappa_{1}}\right\}\le \beta(v)\le 2 \sup\left\{\bar C_1R_\alpha,\bar C_1\left(\frac{7K\lambda_p}{4c_\alpha}\right)^{\frac{1}{k_{\alpha}}}\right\}.  $$
\end{remark}
\medskip

Lastly, we obtain some uniform control on a  continuous set  of homeomorphisms 
$$\tilde \Psi_v(s):=\sum_{i=1}^{N}\Psi_i(sv) v_{i}^2$$ where $v\in U \subset \R^{N,+}$.

Namely, we show that 
\begin{lemma}\label{mutselecedo-lem-diffeo}
Assume that $r_i,(\mu_{ij}), \Psi_i$ satisfies \eqref{mutselecedo-hyp1}, \eqref{mutselecedo-hyp2} and  \eqref{mutselecedo-hyp3}. Assume further that  $\Psi_i(v)=\alpha(v)+\eps \psi_i(v)$ with  $\alpha\in C^1_{loc}$ satisfying \eqref{mutselecedo-hyp1}, \eqref{mutselecedo-hyp2}and $\psi_i \in C^{1}_{loc}$  uniformly bounded. Then there exists $\eps_2$ and $\tau_0>0$ so that for   
 for all  $\eps\le \eps_2$  and for all $\bar v_\eps$ stationary solution of \eqref{mutselecedo-eq-sysgen} we have 
$$ \bar c_3 \le  \tilde \Psi_{\bar v_\eps}(1)\le  \bar C_3 \le 2\bar C_3 \le  \tilde \Psi_{\bar v_\eps}(1+\tau_0).$$
Moreover   there exists $\eps_3$ and $k>0$ so that for all  $\eps\le \eps_3$ we have for all $\bar v_\eps$ stationary solution of \eqref{mutselecedo-eq-sysgen} and  $t,s\in (0,1+\tau_0)$  
$$|t-s|\le k|\tilde \Psi_{\bar v_\eps}(t)-\tilde \Psi_{\bar v_\eps}(s)|.$$
\end{lemma}

\dem{Proof:}

Recall that from the proof of Lemma \ref{mutselecedo-lem-estigen1}  for all $\eps \le \eps_1$, for any  stationary solution of \eqref{mutselecedo-eq-sysgen} $\bar v_\eps$, we  have $\alpha(\bar\o^-)\le \alpha(\bar v_\eps) \le \alpha(\bar\o^+)$. So, for all $\eps\le \eps_1$ and for all stationary solution of \eqref{mutselecedo-eq-sysgen} we have 
$$<\bar v_\eps,\bar v_\eps>(\alpha(\bar \o^-)-\sigma) \le \tilde \Psi_{\bar v_\eps}(1)\le  (\alpha(\bar \o^+)+\sigma)<\bar v_\eps,\bar v_\eps>,$$
where $\sigma=\eps\|\psi\|_{\infty}$. Let us fix $ \eps_2:=\min \{\eps_1,\frac{\lambda_p}{4K\|\psi\|_{\infty}} \}$.  From  \eqref{mutselecedo-eq-esti2-eps} and Lemma \ref{mutselecedo-lem-estigen1} we get 
$$ \bar c_3:=\frac{\bar c_1^2K\lambda_p}{2N}\le \tilde \Psi_{\bar v_\eps}(1)\le  \frac{3\bar C^2_1K\lambda_p}{2}=:\bar C_3,$$
for all $\eps \le\eps_2 $ and for any  stationary solution  of \eqref{mutselecedo-eq-sysgen}.

By using to the monotonicity  of the map $\sigma \to \bar \o^{\pm}$ and Lemma \ref{mutselecedo-lem-estigen1}, we can choose $\eps_2$ smaller if necessary to achieve  for any $\eps\le \eps_2$ and any stationary solution $\bar v_\eps$,  
$$s_0\bar \o_{\eps_1}^-\le s_0\bar v_\eps\le s_0\bar \o^+_{\eps_1}  \quad \text{ for all }\quad s_0>0. $$
The latter inequalities imply that we have for all $\eps \le \eps_2$ and for any stationary solution $\bar v_\eps$,
$$ \tilde \Psi_{\bar v_\eps}(s_0\bar v_\eps)\ge \alpha(s_0\bar \o^-_{\eps_1}) - \frac{\lambda_p}{4K}.$$
Let us fix $s_0$ such that  
$$ \alpha(s_0\bar \o^-_{\eps_1}) - \frac{\lambda_p}{4K}=2\bar C_3.$$
This is always possible since $\alpha$ is monotone increasing and $\lim_{\mu\to \infty} \alpha(\mu\o^-_{\eps_1})=+\infty$.
By construction, $s_0>1$, since  $\tilde \Psi_{\bar v_\eps}$ is monotone increasing and we have 
$$\tilde \Psi_{\bar v_\eps}(s_0)\ge 2\tilde \Psi_{\bar v_\eps}(1).$$
Let us denote  $\tau_0:=s_0-1$.

Now since for each $i$ the function $\Psi_i$ satisfies \eqref{mutselecedo-hyp1}, \eqref{mutselecedo-hyp2} and that $\alpha$ and the $\psi_i$ are $C^1_{loc}(\R^N)$, the function  $\tilde \Psi_{\bar v_\eps}(s):=\sum_{i=1}^N\Psi_i(s\bar v_\eps)\bar v_{\eps,i}^2$ is $C_{loc}^1(\R)$ and monotone increasing. Therefore, for any  stationary solution of \eqref{mutselecedo-eq-sysgen} $\bar v_\eps$, the function $\tilde \Psi_{\bar v_\eps}$ is a $\R^+$ homeomorphism.  

Next, we check that for a fixed $\eps$ and a fixed stationary solution $\bar v_\eps$, the homeomorphism $\tilde \Psi_{\bar v_\eps}$ is a $C^1$ diffeomorphism on $(0,1+\tau_0)\to \tilde \Psi_{\bar v_\eps}((0,1+\tau_0))$. Thanks to the Inverse Function Theorem, to show that  $\tilde \Psi_{\bar v_\eps}$ is a local  $C^1$  diffeomorphism it is sufficient to prove that for all $s \in (0,1+\tau_0)$, $\tilde \Psi_{\bar v_\eps}^{\prime}(s)\neq 0$.
By a straightforward  computation we have:
$$ \left(\Psi_{\bar v_\eps}\right)^{\prime}(s)=<\nabla \alpha(s \bar v_{\eps}),\bar v_\eps>\e(\bar v)+\eps\sum_{i=1}^N<\nabla \psi_i(s\bar v_\eps),\bar v>\bar v_{\eps,i}^2. $$
For all $s\in (0,1+\tau_0)$, and for any  stationary solution $\bar v_\eps \in Q_{(1+\tau_0)\bar C_1}(0)$ we have:
\begin{align*}
&|\nabla \psi_i(s\bar v_\eps)|\le \sup_{ v\in Q_{(1+\tau_0)\bar C_1}(0)}|\nabla \psi_i(v)|:=\zeta_1,\\
&\nabla \alpha(s\bar v_\eps)\le \sup_{v \in Q_{(1+\tau_0)\bar C_1}(0)}\nabla \alpha (v):=\zeta_2 >0,\\
&\nabla \alpha(s\bar v_\eps)\ge \inf_{v\in Q_{(1+\tau_0)\bar C_1}(0)}\nabla \alpha (v):=\zeta_3 >0.
\end{align*}
Therefore, by Lemma \ref{mutselecedo-lem-estigen1} we deduce that 
$$ \left(\Psi_{\bar v_\eps}\right)^{\prime}(s)\ge \frac{\bar c_1^3}{N}\zeta_3- \eps \zeta_1\bar C_1^3 .$$
By choosing $\eps\le \eps_3:=\min \left\{\eps_2,\frac{\bar c_1^3\zeta_3}{2N \zeta_1 \bar C_1^3} \right\}$, we get for all $\eps$, $s\in (0,1+\tau_0)$ and for any   stationary solution of \eqref{mutselecedo-eq-sysgen} $\bar v_\eps$,
\begin{equation*}
0<\frac{\bar c_1^3 \zeta_3}{2N}\le \left(\Psi_{\bar v_\eps}\right)^{\prime}(s)\le \bar C_1^3 \left(\zeta_2+\frac{\zeta_3}{2N}\right)  .
\end{equation*}
From the the latter \textit{ a priori } bounds, we see that for all $\eps\le \eps_3$, $s\in (0,1+\tau_0)$ and for any  stationary solution of \eqref{mutselecedo-eq-sysgen} $\bar v_\eps$, we get the following \textit{a priori} estimate 
\begin{equation*}
  \frac{1}{\bar C_1^3 \left(\zeta_2+\frac{\zeta_3}{2N}\right) }   \le \left(\tilde\Psi_{\bar v_\eps}^{-1}\right)^{\prime}[\tilde\Psi_{\bar v_\eps}(s)]=\frac{1}{\left(\Psi_{\bar v_\eps}\right)^{\prime}(s)}\le \frac{2N}{\bar c_1^2 \zeta_3}.
\end{equation*}

Hence, we deduce that for all $\eps\le \eps_3$, $s,t\in (0,1+\tau_0)$ and for any  stationary solution of \eqref{mutselecedo-eq-sysgen} $\bar v_\eps$, we have 
$$|s-t|\le  k|\tilde \Psi_{\bar v_\eps}(s)-\tilde \Psi_{\bar v_\eps}(t)|,$$
with 
$$k:=\frac{2N}{\bar c_1^2 \zeta_3}.$$
   
\fdem

\subsection{Asymptotic Behaviour}~\\

We are now in position to obtain the asymptotic behaviour of the solution $v_\eps(t)$ as $t$ goes to $+\infty$ for  $\eps \in [0,\eps^*]$, where $\eps^*$ is to be determined  later on. 

Let us fix $\eps \in [0,\eps_1]$ where $\eps_1$ is obtained in Lemma \ref{mutselecedo-lem-estigen1} and let $\bar v_\eps$ be a stationary solution of \eqref{mutselecedo-eq-sysgen}.  For simplicity we denote $<,>$  the standard scalar product in $\R^N$.

As in the proof of Lemma  \ref{mutselecedo-lem-asb}, we start by observing that since $\bar v_\eps \neq 0$  we can write $v_\eps(t):=\lambda(t)\bar v_\eps +h(t)$ with for all $t,$ $<h,\bar v_\eps>=0$. 
In the sequel of this subsection, for more clarity in the presentation we drop the subscript $\eps$ on  $v$ and $\bar v$.

First, we note that from this decomposition we can derive the following equalities:
\begin{align}
&\lambda <\bar v,\bar v>=\sum_{i=1}^Nv_i\bar v_i ,\label{mutselec-eq-sysgen-asb1}\\
&\frac{d\e(h)}{dt}=\frac{d\e(v)}{dt} -2\lambda\lambda' <\bar v,\bar v>\label{mutselec-eq-sysgen-asb2}.
\end{align}
By Lemma  \ref{mutselecedo-lem-estigen1} and \eqref{mutselec-eq-sysgen-asb1}, we can check that  for all $\eps \le \eps_1$ and for all $t>0$, we have 
\begin{equation}\label{mutselecedo-eq-asb-esti1}
\frac{\bar c_2}{\bar C_1^2}\le \lambda(t)\le \frac{N\bar C_2}{\bar c_1^2}.
\end{equation}

Similarly, by using  \eqref{mutselecedo-eq-asb-esti1}, Lemma \ref{mutselecedo-lem-estigen1} and   $\sum_{i=1}^N|h_i|\le \lambda(t)\sum_{i=1}^N\bar v_i +\sum_{i=1}^Nv_i$, we deduce   from $\e(v)=\e(h)+\lambda^2\e(\bar v)$ that 
\begin{align}
& \e(h) \le \frac{N\bar C_2^2}{\bar c_1^2}(N+1),\label{mutselecedo-eq-asb-esti2}\\
 &\sum_{i=1}^N|h_i|\le\frac{N\bar C_2}{\bar c_1}\left(1+\frac{\bar C_1}{\bar c_1}\right).\label{mutselecedo-eq-asb-esti3}
\end{align}

By \eqref{mutselecedo-eq-asb-esti1},\eqref{mutselecedo-eq-asb-esti2},\eqref{mutselecedo-eq-asb-esti3},  Lemma \ref{mutselecedo-lem-estigen1},and by using the Lipschitz regularity of $\Psi_i$ and the Cauchy-Schwartz inequality we can check that for some  constant $C>0$ independent of $\eps$ and $v$
\begin{align}
&\sum _{i=1}^N|\Psi_i(\lambda(t)\bar v)-  \Psi_i(v(t))|\bar v_i^2 \le C\sqrt{\e(h)},\label{mutselecedo-eq-asb-esti4}\\
&\sum_{i=1}^N|(\psi_i(\bar v)-\psi_i(v))\bar v_i h_i|\le C\sqrt{\e(h)}\label{mutselecedo-eq-asb-esti5}.
\end{align}

Indeed  by \eqref{mutselecedo-eq-asb-esti1}--\eqref{mutselecedo-eq-asb-esti3} and Lemma \ref{mutselecedo-lem-estigen1}, $\lambda \bar v$ and $h$ are uniformly bounded and we have
\begin{equation*}
 \sum_{i =1}^{N}|\psi_i(\bar v)-\psi_i(v)|  |h_i|\bar v_i \le\sum_{i =1}^{N}\kappa_i \left(\sum_{j=1}^N |\bar v_j-v_j|\right)  |h_i|\bar v_i,
\end{equation*}
where $\kappa_i$  are the Lipschitz constant of $\psi_i$ on the set $Q_{R_0}(0)$ with $R_0:=\frac{N\bar C_2\bar C_1}{\bar c_1^2}\left(1+2\frac{ \bar C_1  }{ \bar c_1}\right)$. 
From  the decomposition of $v$  and by using the Cauchy-Schwartz inequality, we get 
  \begin{align}
   \sum_{i =1}^{N}| (\psi_i(\bar v)-\psi_i(\lambda(t) \bar v+h)) \bar v_i h_i|&\le \bar \kappa \left( |1-\lambda(t)|\sum_{j=1}^N \bar v_j + \sum_{j=1}^N |h_j| \right)\sum_{i =1}^{N}  |h_i|\bar v_i\\
  &\le N \bar \kappa  \sqrt{\e(h)\e(\bar v)} \left[ |1-\lambda(t)| \sqrt{\e(\bar v)} +\sqrt{\e(h)}\right] \label{mutselecedo-eq-asb-esti6},
  \end{align}
  where $\bar \kappa:=\sup_{i\in \{1,\ldots, N\}}\kappa_i $.  
Therefore, by \eqref{mutselecedo-eq-asb-esti2} and Lemma \ref{mutselecedo-lem-estigen1}, the inequality \eqref{mutselecedo-eq-asb-esti5} holds true for some positive constant $C$ independent of  $\eps \le \eps_1$.  
A similar argumentation holds to get  \eqref{mutselecedo-eq-asb-esti4}.

Next, we show  that 
\begin{lemma} \label{mutselecedo-cla-energy-gen}
There exists $\eps^*\le \min\{\eps_0,\eps_1,\eps_3\},$ so that for all $\eps \le \eps^*$,  $\e(h_{\eps}(t))\to 0$ as $t\to +\infty$.
 \end{lemma}

Assume the lemma holds true, then we can conclude  the proof of Theorem \ref{mutselecedo-thm4} by arguing as follows.

By combining  \eqref{mutselec-eq-sysgen-asb1} and Lemma \ref{mutselecedo-lem-general-id}, we  achieve 
\begin{equation}
\lambda^{\prime}(t)<\bar v, \bar v>=\frac{d}{dt}\sum_{i=1}^Nv_i\bar v_i =\frac{1}{K}\sum_{i=1}^N\Gamma_i v_i\bar v_i.
\end{equation}

 Now by using  $\e(h)\to 0$, we deduce that   $h \to 0$ as $t\to \infty$ and   from the latter equality we are reduced to analyse  the ODE 

$$\lambda'(t)\e(\bar v)=\frac{\lambda(t)}{K}\sum_{i=1}^N(\Psi_i(\bar v)-\Psi_i(\lambda(t) \bar v)) \bar v_i^2+o(1),$$
where 
$$o(1):=\sum _{i=1}^N(\Psi_i(\lambda(t)\bar v)-  \Psi_i(v))\bar v_i v_i,$$     
which by \eqref{mutselecedo-eq-asb-esti4} and  \eqref{mutselecedo-eq-asb-esti5}  satisfies
\begin{equation}
|o(1)|\le \frac{C(1+\eps)}{K}\sqrt{\e(h)}\to 0 \qquad \text{ as }\quad t\to +\infty .\label{mutselec-eq-sysgen-asb-o1}
\end{equation} 

By Lemma  \ref{mutselecedo-lem-estigen1},  for all $\eps \le \eps^*\le \eps_1$ we have  
$\bar c_2\le\sum_{i=1}^N\bar v_iv_i\le \bar C_2$ with  $\bar c_2$ and $\bar C_2$ independent of $\eps$. 
Therefore by \eqref{mutselec-eq-sysgen-asb1} $\frac{\bar c_1}{\e(\bar v)}\le \lambda(t)\le \frac{\bar C_2}{\e(\bar v)}$ for all $t>0$. Thus $\lambda$ satisfies an ODE of the form

$$\lambda'(t)=\frac{\lambda(t)}{K\e(\bar v)}(\theta-\tilde \Psi (\lambda(t)))+\frac{\lambda(t)}{K\e(\bar v)}o(1),$$
which can be rewritten
\begin{equation}
\lambda'(t)=\frac{\lambda(t)}{K\e(\bar v)}(\theta+o(1)-\tilde \Psi (\lambda(t))), \label{mutselec-eq-sysgen-asb-edo}
\end{equation}
where $\theta:=\sum_{i=1}^N\Psi_i(\bar v)\bar v_i^2$ and $\tilde \Psi(s):=\sum_{i=1}^N\Psi_i(s\bar v)\bar v_i^2$. By construction, $\tilde \Psi \in C^{1}_{loc}(\R)$ and  is monotone increasing. 

The above ODE is of logistic type with a perturbation $o(1)\to 0$  with  a non negative initial datum.  Therefore, by a standard argumentation, we see that  $\lambda(t)$ converges to  $\bar \lambda>0$ where $\bar \lambda$ is the unique solution of  $\tilde \Psi(\bar \lambda )=\theta$. By construction, we have $\tilde \Psi(1)=\theta$, so we deduce that $\bar \lambda=1$.
Hence, $\lambda(t)\to 1$ and we can conclude that $v$ converges  to $\bar v$.
\fdem

Let us now turn our attention to the proof of the Lemma \ref{mutselecedo-cla-energy-gen}.

\dem{Proof of Lemma \ref{mutselecedo-cla-energy-gen}:}
First, let us denote $\Gamma_i:=\Psi_i(\bar v)-\Psi_i(v)$.  Then, by
combining  \eqref{mutselec-eq-sysgen-asb1}, \eqref{mutselec-eq-sysgen-asb2}and Lemma \ref{mutselecedo-lem-general-id}, we  achieve 
\begin{equation}
\lambda^{\prime}(t)<\bar v, \bar v>=\frac{d}{dt}\sum_{i=1}^Nv_i\bar v_i =\frac{1}{K}\sum_{i=1}^N\Gamma_i v_i\bar v_i \label{mutselec-eq-sysgen-asb3}
\end{equation}
and 
\begin{align*}
\frac{d\e(h)}{dt}&=-\sum_{i,j =1}^{N}\mu_{ij} \bar v_i \bar v_j\left( \frac{h_j}{\bar v_j}-\frac{h_j}{\bar v_j}\right)^2+\frac{2}{K}\sum_{i =1}^{N}\Gamma_i h_i^2 +\frac{2\lambda}{K}\sum_{i =1}^{N}\Gamma_i h_i\bar v_i,\\
&=-\sum_{i,j =1}^{N}\mu_{ij} \bar v_i \bar v_j\left( \frac{h_j}{\bar v_j}-\frac{h_j}{\bar v_j}\right)^2+\frac{2}{K}\sum_{i =1}^{N}\Gamma_i v_i h_i.
\end{align*}
Therefore using the definition of $\Psi_i$ and with the notation $\gamma_i:=\psi_i(\bar v) -\psi_i(v)$, we have 
\begin{multline*}
\frac{d\e(h)}{dt}=-\sum_{i,j =1}^{N}\mu_{ij} \bar v_i \bar v_j\left( \frac{h_j}{\bar v_j}-\frac{h_j}{\bar v_j}\right)^2+\frac{2}{K}(\alpha(\bar v)-\alpha(v))\e(h)\\ +\frac{2\eps}{K}\sum_{i =1}^{N}\gamma_i  h_i^2+\frac{2\eps\lambda(t)}{K}\sum_{i =1}^{N}\gamma_i  \bar v_i h_i,
\end{multline*}
which implies that 
\begin{multline}
\frac{d\e(h)}{dt}\le-\sum_{i,j =1}^{N}\mu_{ij} \bar v_i \bar v_j\left( \frac{h_j}{\bar v_j}-\frac{h_j}{\bar v_j}\right)^2+\frac{2}{K}(\alpha(\bar v)+2\eps\|\psi\|_{\infty}-\alpha(v))\e(h)\\ +\frac{2\eps\lambda(t)}{K}\sum_{i =1}^{N}\gamma_i  \bar v_i h_i.\label{mutselec-eq-sysgen-asb4}
\end{multline}

By \eqref{mutselec-eq-sysgen-asb3}, using the definition of $\Psi_i$ we also have 
\begin{align*}
\frac{d}{dt}\sum_{i=1}^N\bar v_i v_i&=\frac{1}{K}(\alpha(\bar v) -\alpha(v))\sum_{i=1}^N\bar v_i v_i +\frac{\eps}{K}\sum_{i=1}^N\gamma_i\bar v_i v_i,\\
&\ge \frac{1}{K}(\alpha(\bar v)-2\eps\|\psi\|_{\infty} -\alpha(v))\sum_{i=1}^N\bar v_i v_i. 
\end{align*}
By Lemma  \ref{mutselecedo-lem-estigen1},  for all $\eps \le \eps_1,$  
$\bar c_2\le\beta(v)=\sum_{i=1}^N\bar v_iv_i\le \bar C_2$ with  $\bar c_2$ and $\bar C_2$ independent of $\eps$. So we have
$$\frac{d}{dt}\log(\beta(v(t))) \ge  \frac{1}{K}(\alpha(\bar v)-2\eps\|\psi\|_{\infty} -\alpha(v)),$$ 
which combined with \eqref{mutselec-eq-sysgen-asb4} leads to 
  \begin{equation*}
\frac{d\e(h)}{dt}\le-\sum_{i,j =1}^{N}\mu_{ij} \bar v_i \bar v_j\left( \frac{h_j}{\bar v_j}-\frac{h_j}{\bar v_j}\right)^2+\frac{d}{dt}\log(\beta^2(v(t))) \e(h)+\frac{8\eps\|\psi\|_{\infty}}{K}\e(h) +\frac{2\eps\lambda(t)}{K}\sum_{i =1}^{N}\gamma_i  \bar v_i h_i.
  \end{equation*}
By using the functional inequality, Lemma \ref{mutselecedo-lem-fcineq}, and rearranging the terms in the above inequality we get
 \begin{equation}
\frac{d\e(h)}{dt}-\e(h)\frac{d}{dt}\log(\beta^2(v(t)))  \le\left(-C_1(\bar v_\eps)+\frac{8\eps\|\psi\|_{\infty}}{K}\right)\e(h) +\frac{2\eps\lambda(t)}{K}\sum_{i =1}^{N}\gamma_i  \bar v_i h_i,\label{mutselec-eq-sysgen-asb6}
  \end{equation}
where $C_1(\eps)$ is the second largest eigenvalue of some associated linear problem.

Now, we estimate the last term of the above inequality. Recall that by  \eqref{mutselecedo-eq-asb-esti6} we have 
 
$$ 
 \sum_{i =1}^{N}|\gamma_i  \bar v_i h_i|\le N \bar \kappa  \sqrt{\e(h)\e(\bar v)} \left[ |1-\lambda(t)| \sqrt{\e(\bar v)} +\sqrt{\e(h)}\right].$$

  By combining the above estimate and  \eqref{mutselec-eq-sysgen-asb6}, we achieve
  \begin{equation}
\frac{d\e(h)}{dt}-\e(h)\frac{d}{dt}\log(\beta^2(v(t)))  \le\left(-C_1(\bar v_\eps)+\eps C_5 \right)\e(h) +\eps C_4|1-\lambda(t)|\sqrt{\e(h)},  \label{mutselec-eq-sysgen-asb7}
  \end{equation}
  
   where $C_5:=\frac{2}{K}\left(4\|\psi\|_{\infty}+\frac{N^2\bar \kappa \bar C_2}{\bar c_1}\right)$ and $C_4:=\frac{2\bar C_2N\bar\kappa}{K}$.

Recall that $C_1(\bar v)$ is an eigenvalue, therefore for $\eps$ small enough, say $\eps \le \eps_4$,  since $ \bar v\in  Q_{C\sqrt{\eps}}(\bar \o^+)$  and  $\bar \o^+\to \bar \o$ as $\eps\to 0$ one has 
$C_1(\bar v)\ge \frac{C(\o^+)}{2}$.

With the latter estimate and by choosing  $\eps\le \min( \eps_1,\eps_4)$ smaller if necessary, we get

\begin{equation}
\frac{d\e(h)}{dt}-\e(h)\frac{d}{dt}\log(\beta^2(v(t)))  \le -\frac{C_1(\bar \o^+)}{4}\e(h) +\eps C_4|1-\lambda(t)|\sqrt{\e(h)}.  \label{mutselec-eq-sysgen-asb8}
  \end{equation}

The proof now will follow several steps:
\subsubsection*{Step One:}
Since by \eqref{mutselecedo-eq-asb-esti1} we have $|1-\lambda(t)|\le \left[1+\frac{N\bar C_2}{\bar c_1^2}\right]$ for all $t$,
we  claim that 
\begin{claim}
For all $\eps \le \min( \eps_1,\eps_4) $, there exists $t_0$ so that for all $t\ge t_0$ we have 
$$\sqrt{\e(h)}\le 2\eps \left(\frac{\bar C_2}{\bar c_2}\right)\frac{4C_4\left[1+\frac{N\bar C_2}{\bar c_1^2}\right]}{C_1(\bar \o^+)}$$
\end{claim}

Indeed, by \eqref{mutselec-eq-sysgen-asb1}, and Lemma \ref {mutselecedo-lem-estigen1} we have 
\begin{equation}
\frac{d\e(h)}{dt}-\e(h)\frac{d}{dt}\log(\beta^2(v(t)))  \le -\frac{C_1(\bar \o^+)}{4}\e(h) +\eps C_4\left[1+\frac{N \bar C_2}{\bar c_1^2}\right]\sqrt{\e(h)}  \label{mutselec-eq-sysgen-asb9}
  \end{equation}

 and we can check that there exists $t_0>0$ so that $ \sqrt{\e(h(t_0))}\le 2\frac{\eps 4C_4\left[1+\frac{N \bar C_2}{\bar c_1^2}\right]}{C_1(\bar \o^+)}$.  If not, then  for all $t>0$ $ \sqrt{\e(h(t))}>2\frac{\eps 4C_4\left[1+\frac{N \bar C_2}{\bar c_1^2}\right]}{C_1(\bar \o^+)}$ and by dividing  \eqref{mutselec-eq-sysgen-asb9} by $\sqrt{\e(h)}$  and 
 rearranging the terms, we get that 
\begin{equation}
\sqrt{\e(h)}\frac{d}{dt}\log\left(\frac{\e(h)}{\beta^2(v(t))}\right)  \le -\frac{C_1(\bar \o^+)}{4}\sqrt{\e(h)} +\eps C_4\left[1+\frac{N \bar C_2}{\bar c_1^2}\right]< -\eps C_4\left[1+\frac{N \bar C_2}{\bar c_1^2}\right].  \label{mutselec-eq-sysgen-asb10}
  \end{equation}
  Thus $ F(t):=\log\left(\frac{\e(h)}{\beta^2(v(t))}\right)$ is a decreasing function which by Lemma \ref{mutselecedo-lem-estigen1}, is bounded from below.  Moreover, we have for all $t>0$ $ \sqrt{\e(h(t))}>2\frac{\eps 4C_4\left[1+\frac{N \bar C_2}{\bar c_1^2}\right]}{C_1(\bar \o^+)}$. Therefore $F$ converges as $t$ tends to $+\infty$ and $\frac{dF}{dt} \to 0$. 
  Thus for $t$ large enough, we get the contradiction 
  $$-\frac{1}{2}\eps C_4\left[1+\frac{N \bar C_2}{\bar c_1^2}\right]\le \sqrt{\e(h)}\frac{d}{dt}\log\left(\frac{\e(h)}{\beta^2(v(t))}\right)  \le \eps C_4\left[1+\frac{N \bar C_2}{\bar c_1^2}\right].$$

Let $\Sigma$ be  the set  $\Sigma:=\left\{t>t_0 | \sqrt{\e(h(t))}>  2\frac{\eps 4C_4\left[1+\frac{N \bar C_2}{\bar c_1^2}\right]}{C_1(\bar \o^+)}\right\}$.  
Assume that $\Sigma$ is non empty otherwise the claim is proved  since $\frac{\bar C_2}{\bar c_2}>1$.  Let us  denote   $t^*:=\inf\Sigma $. So at  $t^*$ we have $\sqrt{\e(h(t^*))} = 2\frac{\eps 4C_4\left[1+\frac{N \bar C_2}{\bar c_1^2}\right]}{C_1(\bar \o^+)}$. 
Again,  by dividing\eqref{mutselec-eq-sysgen-asb9}  by $\sqrt{\e(h)}$  and rearranging the terms,   on the set $\Sigma$  we have
\begin{equation}
\sqrt{\e(h)}\frac{d}{dt}\log\left(\frac{\e(h)}{\beta^2(v(t))}\right)  \le -\frac{C_1(\bar \o^+)}{4}\sqrt{\e(h)} +\eps C_4\left[1+\frac{N \bar C_2}{\bar c_1^2}\right]\le 0.  \label{mutselec-eq-sysgen-asb11}
  \end{equation}

Thus $\log\left(\frac{\e(h)}{\beta^2(v(t))}\right)$ is a decreasing function of $t$ for all $t \in \Sigma$. By arguing on each connected component of $\Sigma$, we can check that  for all $t \in \Sigma$
$$ \sqrt{\e(h(t))}\le \frac{\beta(v(t)) }{\bar c_2}\frac{\eps 4C_4\left[1+\frac{N \bar C_2}{\bar c_1^2}\right]}{C_1(\bar \o^+)}.$$
Therefore by using Lemma \ref{mutselecedo-lem-estigen1} for $t\ge t^*$ we have
$$\sqrt{\e(h)}\le 2\frac{\bar C_2}{\bar c_2}\frac{\eps 4C_4\left[1+\frac{N \bar C_2}{\bar c_1^2}\right]}{C_1(\bar \o^+)}.$$

Hence, since $\frac{\bar C_2}{\bar c_2}>1$ we get for all $t\ge t_0$, 
$$\sqrt{\e(h)}\le 2\frac{\bar C_2}{\bar c_2}\frac{\eps 4C_4\left[1+\frac{N \bar C_2}{\bar c_1^2}\right]}{C_1(\bar \o^+)}.$$
\fdem

\subsubsection*{Step Two:}

First, we define some constant quantities: 
\begin{align}
&\delta_0:=2\frac{\bar C_2}{\bar c_2}\frac{4C_4\left[1+\frac{N \bar C_2}{\bar c_1^2}\right]}{C_1(\bar \o^+)},\label{mutselecedo-eq-sysgen-asb12}\\
&d:=\left[1+\frac{N \bar C_2}{\bar c_1^2}\right],\label{mutselecedo-eq-sysgen-asb13}\\
&\eps^*:=\min\left\{\eps_1,\eps_2, \eps_3,\eps_4,\frac{d}{4kC\delta_0},\frac{\bar c_3}{4C\delta_0}\right\},
\end{align}
where the constants $C,k$ are respectively   defined in \eqref{mutselecedo-eq-asb-esti4},\eqref{mutselecedo-eq-asb-esti5} and in the  Lemma \ref{mutselecedo-lem-diffeo}.

By the previous step, we have that for all $\eps\le \eps^*, $ all   $t\ge t_0$
$$ \sqrt{\e(h)}\le \eps\delta_0.$$
We claim that 
\begin{claim}\label{mutselecedo-cla-steptwo}
For all $\eps\le \eps^*$   there exists $t_ {\eps\delta_0}$ such that for all $t\ge t_{\eps\delta_0}$
$$\sqrt{\e(h(t))}\le \frac{\eps\delta_0}{2}.$$
 \end{claim}
 \dem{Proof:}	
First, we  can check that for $\eps\le\eps^*$  there exists $t^*$  so that for all $t\ge t^*$
$$|1-\lambda(t)|\le 2k C\eps\delta_0.$$

Indeed, by \eqref{mutselec-eq-sysgen-asb-edo} and \eqref{mutselec-eq-sysgen-asb-o1},  we have
\begin{equation}
\lambda'(t)=\frac{\lambda(t)}{K\e(\bar v)}(\theta+o(1)-\tilde \Psi (\lambda(t))),
\end{equation}  
 with $ o(1)\le C\eps\delta_0$ for $t \ge t_0$.
 Let $\lambda_{\pm C\eps\delta_0}$ be the solution of the ODE
 \begin{equation}
\lambda_{\pm C\eps\delta_0}'(t)=\frac{\lambda_{\pm C\eps\delta_0}(t)}{K\e(\bar v)}( \tilde \Psi_{\bar v_\eps} (1)\pm C\eps\delta_0-\tilde \Psi_{\bar v_\eps} (\lambda_{\pm C\eps\delta_0}(t))).
\end{equation}  
Since $\eps\le \frac{\bar c_3}{4C\delta_0}$, for $t\ge t_0$ we have 
$$\frac{3\bar c_3}{4}\le \tilde \Psi_{\bar v_\eps} (1)\pm C\eps\delta_0 \le \frac{5 \bar C_3}{4}. $$  
Therefore $\lambda_{\pm C\eps\delta_0}\to \bar\lambda_{\pm C\eps\delta_0} $ where $\bar\lambda_{\pm C\eps\delta_0}$ are the unique positive solutions of $\tilde \Psi_{\bar v_\eps}( \bar\lambda_{\pm C\eps\delta_0})=\tilde \Psi_{\bar v_\eps} (1) \pm C\eps\delta_0$.

Thanks to the strict  monotonicity of $\tilde \Psi_{\bar v_\eps} $, we also have 
 $$\bar \lambda_{-2C\eps\delta_0}< \bar \lambda_{-C\eps\delta_0}< \bar \lambda_{C\eps\delta_0}< \bar \lambda_{2C\eps\delta_0}, $$ 
where $\bar\lambda_{\pm 2C\eps\delta_0}$ are the unique positive solutions of $\tilde \Psi_{\bar v_\eps}( \bar\lambda_{\pm 2C\eps\delta_0})=\tilde \Psi_{\bar v_\eps} (1) \pm 2C\eps\delta_0$.

Since  $\eps\le \eps^*\le \frac{\bar c_3}{4C\delta_0}$, it follows that
$$ \frac{\bar c_3}{2}\le\tilde \Psi_{\bar v_\eps} (1) \pm 2C\eps\delta_0\le \frac{3 \bar C_3}{2}$$
 and therefore, by Lemma \eqref{mutselecedo-lem-diffeo},  we have 
\begin{equation}
 0<\bar \lambda_{\pm2C\eps\delta_0}<1+\tau_0. \label{mutselecedo-eq-tau0}
 \end{equation}

Now, recall that   for $t\ge t_0$, $\lambda(t)$ satisfies 
\begin{align*}
&\lambda'(t)\le \frac{\lambda(t)}{K\e(\bar v)}(\tilde \Psi_{\bar v_\eps}(1) +C\eps\delta_0 -\tilde \Psi_{\bar v_\eps} (\lambda(t))),\\
&\lambda'(t)\ge \frac{\lambda(t)}{K\e(\bar v)}(\tilde \Psi_{\bar v_\eps}(1) -C\eps\delta_0 -\tilde \Psi_{\bar v_\eps} (\lambda(t))).
\end{align*} 
 Thus, by a standard  argumentation, we can  show  that for $t\ge t^*$ we have 
$$\bar \lambda_{-2C\eps\delta_0}\le \lambda(t) \le \bar \lambda_{+2C\eps\delta_0}.$$

Therefore, for $t\ge t^*$ we have 
$$|1-\lambda(t)|\le \sup \{|1-\bar \lambda_{-2C\eps\delta_0}|,|1-\lambda_{+2C\eps\delta_0}| \}.$$
 By \eqref{mutselecedo-eq-tau0} and  Lemma \ref{mutselecedo-lem-diffeo},  since $\eps\le \eps^*\le\eps_3$,    we deduce that for $t\ge t^*$, 
$$
 |1-\lambda(t)|\le k\sup \{|\tilde \Psi_{\bar v_\eps}(1)-\tilde \Psi_{\bar v_\eps}(\bar \lambda_{-2C\eps\delta_0})|,|\tilde \Psi_{\bar v_\eps}(1)-\tilde \Psi_{\bar v_\eps}(\lambda_{+2C\eps\delta_0})| \} \le 2kC\eps\delta_0.
 $$

From the latter estimate,  by using \eqref{mutselec-eq-sysgen-asb8}, we see that for $t\ge t^*$

\begin{equation}
\frac{d\e(h)}{dt}-\e(h)\frac{d}{dt}\log(\beta^2(v(t)))  \le -\frac{C_1(\bar \o^+)}{4}\e(h) +2\eps^2 C_4kC\delta_0\sqrt{\e(h)}.  \label{mutselec-eq-sysgen-asb14}
  \end{equation}
  By following the argumentation of Step one,  we can show that there exists $ t_{\eps\delta_0}$ such that for   $t\ge t_{\eps\delta_0}$ we have 
$$\sqrt{\e(h)}\le \frac{2kC\eps^2\delta_0^2}{d}.$$ 

Hence, we have for all $t\ge t_{\eps\delta_0}$  
$$\sqrt{\e(h)}\le \frac{\eps\delta_0}{2},$$
since   $\eps\delta_0<\frac{d}{4kC}$.
\fdem

\subsubsection*{Step Three:}
Since for all $t\ge t_{\eps\delta_0}$, 
$$\sqrt{\e(h(t))}\le \frac{\eps\delta_0}{2},$$
 by arguing as in the proof of Claim \ref{mutselecedo-cla-steptwo}, we show that for all $\eps \le \eps^*$ there exists $t_{\eps\frac{\delta_0}{2}}$ so that for all $t\ge t_{\eps\frac{\delta_0}{2}}$
$$\sqrt{\e(h)}\le \frac{\eps\delta_0}{4}.$$ 

By reproducing inductively the above argumentation,  for all $\eps \le \eps^*$ we can construct a sequence  $(t_{n})_{n\in \N}$ so that 
for all $t\ge t_n$ we have 
$$ \sqrt{\e(h(t))}\le \frac{\eps\delta_0}{2^n}.$$
Hence, for all $\eps \le \eps^*$ we deduce that  
$$\lim_{t\to \infty}\e(h(t)) \to 0.$$

\fdem

\bigskip

\noindent \textbf{Acknowledgements.}  J. Coville thanks the member of the  project ERBACE of INRIA for early discussion on this subject.  F. Fabre is supported by the project VirAphid of the \textit{Agence Nationale de la Recherche} (project ANR-10-STRA-0001).

\bibliographystyle{plain}
\bibliography{mutselecedo.bib}

\end{document}